\documentclass[11pt,letterpaper]{amsart}
\pdfoutput=1  
\usepackage{eucal}
\usepackage{graphicx}
\usepackage{amsopn,amssymb}
\usepackage{xypic}
\xyoption{all}

\newcommand{\dzbardz}{\frac{d\overline{z}}{dz}}
\newcommand{\slashdzbardz}{d\overline{z}/dz}
\newcommand{\dbard}[1]{\frac{d\overline{#1}}{d{#1}}}
\newcommand{\slashdbard}[1]{d\overline{#1}/d{#1}}
\newcommand{\dtzero}{\left.\frac{d}{dt}\right|_{t=0^+}}
\newcommand{\dt}{\frac{d}{dt}}

\DeclareMathOperator{\Span}{\mathrm{span}}
\newcommand{\param}{\centerdot}  
\newcommand{\sboldpoint}[1]{\smallskip\par\noindent\textbf{#1}}

\newcommand{\id}{\mathrm{Id}}

\newcommand{\R}{\mathbb{R}}
\newcommand{\C}{\mathbb{C}}
\newcommand{\N}{\mathbb{N}}
\newcommand{\Q}{\mathcal{Q}}
\newcommand{\sC}{\mathcal{C}}
\newcommand{\sS}{\mathcal{S}}
\newcommand{\sH}{\mathcal{H}}
\newcommand{\sM}{\mathcal{M}}

\DeclareMathOperator{\im}{Im}
\DeclareMathOperator{\re}{Re}
\renewcommand{\Re}{\re}
\renewcommand{\Im}{\im}

\newcommand{\T}{\mathcal{T}}
\newcommand{\B}{\mathcal{B}}
\renewcommand{\P}{\mathcal{P}}
\newcommand{\QF}{\mathcal{QF}}
\newcommand{\F}{\mathcal{F}}
\newcommand{\CP}{\mathbb{CP}}
\renewcommand{\H}{\mathbb{H}}
\newcommand{\ML}{\mathcal{ML}}
\newcommand{\GL}{\mathcal{GL}}

\newcommand{\PML}{\mathbb{P} \! \mathcal{M\!L}}

\renewcommand{\Tilde}[1]{\widetilde{#1}}

\DeclareMathOperator{\PSL}{\mathrm{PSL}}
\DeclareMathOperator{\sh}{sh}
\DeclareMathOperator{\gr}{gr}
\DeclareMathOperator{\pr}{pr}
\DeclareMathOperator{\Gr}{Gr}
\DeclareMathOperator{\area}{Area}

\DeclareMathOperator{\Mod}{Mod}
\DeclareMathOperator{\supp}{supp}

\numberwithin{equation}{section}

\theoremstyle{plain}
\newtheorem{thm}{Theorem}[section]
\newtheorem{cor}[thm]{Corollary}

\newtheorem{lem}[thm]{Lemma}

\newtheorem*{main_corollary}{Corollary \ref{cor:graph}}

\newtheorem{example}{Example}

\theoremstyle{definition}
\newtheorem*{question}{Question}

\theoremstyle{remark}
\newtheorem*{remark}{Remark}
\newtheorem*{remarksenv}{Remarks}

\newenvironment{remarks}{\begin{remarksenv}\indent\par\begin{enumerate}}{\end{enumerate}\end{remarksenv}}
\newenvironment{rmenumerate}{\begin{enumerate}
}{\end{enumerate}}

\begin{document}

\title{Projective structures, grafting, and measured laminations}

   \author[D. Dumas]{David Dumas}
   \address{Brown University}
   \thanks{The first author was partially supported by
   an NSF postdoctoral research fellowship.}
   \email{ddumas@math.brown.edu}
   \urladdr{http://www.math.brown.edu/\textasciitilde ddumas/}
   \author[M. Wolf]{Michael Wolf}
   \address{Rice University}
   \thanks{The second author was partially supported by NSF grants 
DMS-0139877 and DMS-0505603.}
   \email{mwolf@math.rice.edu}
   \urladdr{http://www.math.rice.edu/\textasciitilde mwolf/}

   \date{November 30, 2007}

\begin{abstract}
We show that grafting any fixed hyperbolic surface defines a
homeomorphism from the space of measured laminations to Teichm\"uller
space, complementing a result of Scannell-Wolf on grafting by a fixed
lamination.  This result is used to study the relationship between the
complex-analytic and geometric coordinate systems for the space of
complex projective ($\CP^1$) structures on a surface.

We also study the rays in Teichm\"uller space associated to the grafting
coordinates, obtaining estimates for extremal and hyperbolic length
functions and their derivatives along these grafting rays.
\end{abstract}

\maketitle

\section{Introduction}

In this paper we compare two perspectives on the theory of complex
projective structures on surfaces by studying the grafting map of a
hyperbolic surface.

A complex projective (or $\CP^1$) structure on a compact surface $S$
is an atlas of charts with values in $\CP^1$ and M\"obius transition
functions.  Let $\P(S)$ denote the space of (isotopy classes of)
marked complex projective structures on $S$, and let $\T(S)$ be the
Teichm\"uller space of (isotopy classes of) marked complex structures
on $S$.  Because M\"obius maps are holomorphic, there is a forgetful
projection $\pi:\P(S)\to\T(S)$.

An analytic tradition, having much in common with univalent function
theory, parameterizes the fiber $\pi^{-1}(X)$ using the Schwarzian
derivative, identifying $\P(S)$ with the total space of the bundle
$\Q(S)\to\T(S)$ of holomorphic quadratic differentials.

A second, more synthetic geometric description of $\P(S)$ is due to
Thurston, and proceeds through the 
operation of \emph{grafting} -- a construction which 
traces its roots back at least to
Klein \cite[\S50, p. 230]{Kl33}, with a modern history developed by
many authors (\cite{Ma69}, \cite{He75}, \cite{ST83}, \cite{Gol87},
\cite{GKM00},\cite{tanigawa:grafting}, 
\cite{mcmullen:earthquakes}, \cite{scannell-wolf:grafting}).
The simplest example of grafting may be described as follows.

Start with a hyperbolic surface $X \in \T(S)$ and a simple closed
geodesic $\gamma$ on $X$; then construct a new surface by removing
$\gamma$ from $X$ and replacing it with the Euclidean cylinder $\gamma
\times [0,t]$.  The result is $\Gr_{t \gamma}X$, the \emph{grafting of
$X$ by $t \gamma$}, which is a surface with a ($C^{1,1}$ Riemannian)
metric composed of alternately flat or hyperbolic pieces.
Furthermore, $\Gr_{t \gamma} X$ has a canonical projective structure
that combines the Fuchsian uniformization of $X$ and the Euclidean
structure of the cylinder $\gamma \times [0,t]$ (for details, see
\cite[\S1]{scannell-wolf:grafting} \cite[\S2]{tanigawa:grafting}
\cite{kamishima-tan:grafting}).

Thurston showed that grafting extends naturally
from weighted simple closed geodesics to the space $\ML(S)$ of
measured geodesic laminations, and thus defines a map
\begin{equation}\label{eqn:Thurstonhomeo}
\Gr: \ML(S) \times \T(S) \to \P(S).
\end{equation}
Moreover, this map is a homeomorphism; for a proof of this result, see
\cite{kamishima-tan:grafting}.

A natural problem is to relate the analytic and geometric perspectives
on the space of projective structures, for example by comparing the
product structure of $\ML(S) \times \T(S) \simeq \P(S)$ to the bundle
structure induced by the projection $\pi : \P(S) \to \T(S)$.

\subsection*{Results on grafting.} 
We compare these two perspectives by studying the \emph{conformal
grafting map} $\gr = \pi \circ \Gr : \ML(S)\times\T(S)\to\T(S)$,
i.e. $\gr_{\lambda}X$ is the conformal structure which underlies the
projective structure $\Gr_{\lambda} X$.  Fixing either of the two
coordinates we have the $X$-grafting map $\gr_{\param} X : \ML(S) \to
\T(S)$ and the $\lambda$-grafting map $\gr_\lambda : \T(S) \to \T(S)$.
These maps reflect how the base coordinate of the complex-analytic
fibration $\pi : \P(S) \to \T(S)$ is related to the geometric product
structure $\ML(S) \times \T(S) \simeq \P(S)$.  Our main result is

\begin{thm} \label{thm:main}
For each $X\in\T(S)$, the $X$-grafting map
$\gr_{\param}X:\ML(S)\to\T(S)$ is a bitangentiable homeomorphism.
\end{thm}

Momentarily deferring a brief discussion of the term
\emph{bitangentiable homeomorphism}, we note that this theorem is a
natural complement to the result of Scannell-Wolf on the
$\lambda$-grafting map.

\begin{thm}[{Scannell-Wolf \cite[Thm.~A]{scannell-wolf:grafting}}]
\label{thm:scannell-wolf}
For each $\lambda \in \ML(S)$, the map $\gr_\lambda : \T(S) \to \T(S)$
is a real-analytic diffeomorphism.
\end{thm}

The discrepancy between diffeomorphism and bitangentiable homeomorphism
in Theorems \ref{thm:main} and \ref{thm:scannell-wolf} is related to
the lack of a natural differentiable structure on $\ML(S)$.  Bonahon
showed that grafting is differentiable in the weak sense of being
\emph{tangentiable}; see \S2 below or \cite{bonahon:variations} for
details.

Returning to the original problem of comparing different coordinate
systems for $\P(S)$, Theorems \ref{thm:main} and
\ref{thm:scannell-wolf} can be used to study the fiber $P(X) =
\pi^{-1}(X)$ and its relation to the grafting coordinates.  Let us
denote the two factors of the map $\Gr^{-1} : \P(S) \to \ML(S) \times
\T(S)$ by
\begin{equation*}
\Gr^{-1}(Z) = ( p_{\ML}(Z), p_{\T}(Z) ).
\end{equation*}
Thus the maps $p_{\ML} : \P(S) \to \ML(S)$ and $p_{\T} : \P(S) \to
\T(S)$ send a projective structure to one of its two grafting
coordinates, and we think of them as projections.  We prove:

\begin{cor}
\label{cor:graph} For each $X\in\T(S)$ the restriction  
$\left.p_{\ML}\right|_{P(X)} : P(X)\to\ML(S)$ is a bitangentiable
  homeomorphism, and $\left.p_{\T}\right|_{P(X)} : P(X) \to \T(S)$ is a $C^1$
  diffeomorphism.
\end{cor}

\begin{figure}
\begin{center}
\includegraphics{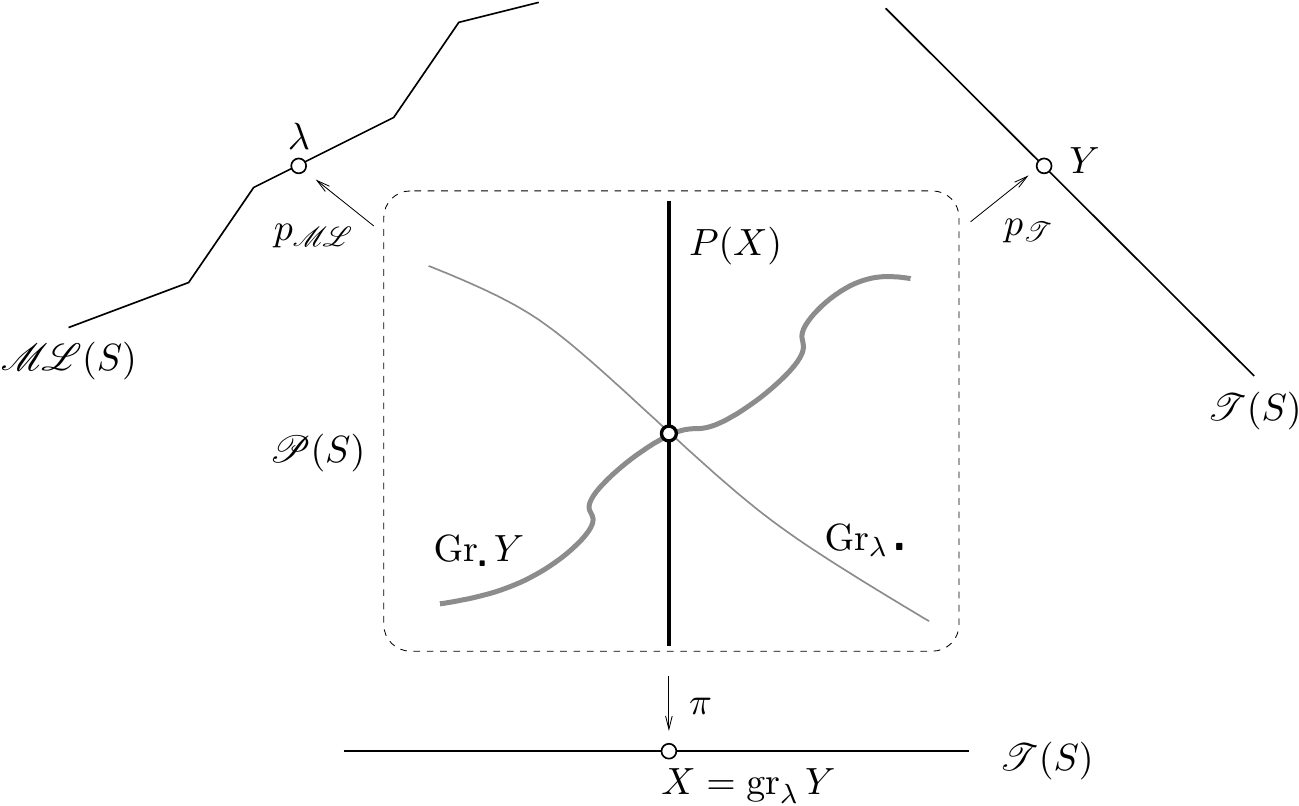}
\caption{The bundle $\pi : \P(S) \to \T(S)$ of $\CP^1$ structures over
  Teichm\"uller space and the product structure $\Gr : \ML(S) \times
  \T(S) \simeq \P(S)$ induced by grafting.}
\label{fig:projections}
\end{center}
\end{figure}

This corollary improves the existing regularity results for these
projection maps, from which it was known that that
$\left.p_{\ML}\right|_{P(X)}$ is a homeomorphism (a corollary of
Theorem \ref{thm:scannell-wolf}, see \cite[\S4]{dumas:schwarzian}) and
that $\left.p_{\T}\right|_{P(X)}$ is a proper $C^1$ map of degree $1$
(see \cite[Thm.~3]{bonahon:variations} and \cite[Lem.~7.6,
Thm.~1.1]{dumas:antipodal}).  The relationship between the maps
$p_{\T}$, $p_{\ML}$, and $\pi$ is represented schematically in Figure
\ref{fig:projections}.

One can also apply Theorem \ref{thm:main} and Corollary
\ref{cor:graph} to study the \emph{pruning map} (the inverse of
grafting)  and the parameterization of quasi-Fuchsian manifolds by
their convex hull geometry.  We explore these directions in
\S\ref{sec:applications}.

\subsection*{Methods.}
From a general perspective, the proof of the main theorem relies on
relating two established techniques in hyperbolic geometry.  The first
is the analytic study of the prescribed curvature (Liouville) equation
(see e.g. \cite{wolpert:WPcurvature}
\cite{takhtajan-zograf:liouville}, \cite{tromba:teichmullerbook}), and
the second is the complex duality between bending and twisting (see
e.g. \cite{wolpert:FNtwistformula}, \cite{wolpert:FNdeformation},
\cite{wolpert:symplectic}, \cite{platis}, \cite{series:Kminima}, \cite{series:smallbending}) and its generalization to complex earthquakes
as a duality between grafting and shearing (see
\cite{bonahon:variations}, \cite{mcmullen:earthquakes},
\cite{epstein-marden-markovic:quasiconformal}).

The first strand occurs in the proof in
\cite{scannell-wolf:grafting} of Theorem~\ref{thm:scannell-wolf},
where standard geometric analytic techniques were applied to the
curvature equation to understand how the $\lambda$-grafting map
changed under small perturbations of the Riemann surface.  Such
techniques might be applicable to the analogous problem (of the main
Theorem~\ref{thm:main}) of understanding how the $X$-grafting map
varies under small perturbations of the measured lamination, but it
would necessarily be more involved, due to the local structure of the
space $\ML(S)$ being more complicated than that of the space $\T(S)$.

Fortunately, most of the required details for this study of $\ML(S)$
are already in the literature: here we make heavy use of Bonahon's
work (following Thurston \cite{thurston:stretch}) on the deformation
theory of $\ML(S)$ (see \cite{bonahon:shearing}
\cite{bonahon:geodesic} \cite{bonahon:transverse}
\cite{bonahon:variations}). In particular, the crux of our proof
relies on Bonahon's observation that there is a sense in which
infinitesimal grafting is complex linear. This complex analyticity, in
keeping with the second tradition discussed above, implies that a
study of the effect on grafting of infinitesimally changing the
measured lamination is, by duality, a study of the effect on grafting
of infinitesimally shearing the hyperbolic surface.  Thus we may apply
the analysis in the proof of Theorem~\ref{thm:scannell-wolf} to the
problem of the main Theorem~\ref{thm:main}.

\subsection*{Grafting rays.}  In a final section,
we study the coordinate system on Teichm\"uller space induced by the
$X$-grafting homeomorphism $\ML(S) \to \T(S)$, and analyze the
behavior of extremal and hyperbolic length functions on \emph{grafting
rays}--paths in $\T(S)$ of the form $t \mapsto \gr_{t \lambda} X$.  We
show:

\begin{thm}
\label{thm:extremal-length-intro}
For each $X \in \T(S)$ and $\lambda \in \ML(S)$ the extremal length of
$\lambda$ on $\gr_{t \lambda}X$ is monotone decreasing for all $t \gg
0$ and is asymptotic to $\frac{\ell(\lambda,X)}{t}$, where
$\ell(\lambda,X)$ is the hyperbolic length of $\lambda$ on $X$.
\end{thm}

\begin{thm}
\label{thm:hyp-length-intro}
For each $X \in \T(S)$ and any simple closed hyperbolic geodesic
$\gamma \in \ML(S)$, the hyperbolic length of $\gamma$ on $\gr_{t
\gamma} X$ is monotone decreasing for all $t \gg 0$ and is asymptotic
to $\pi \frac{\ell(\lambda,X)}{t}$.
\end{thm}

The monotonicity and asymptotic behavior described in Theorems
\ref{thm:extremal-length-intro} and \ref{thm:hyp-length-intro} are
combined with explicit estimates on the derivatives of length
functions in Theorems \ref{thm:extremal-length} and
\ref{thm:hyp-length} below.

\subsection*{Organization of the paper.\\}
\sboldpoint{\S 2} presents the infinitesimal version of the main
theorem (Theorem \ref{thm:infinitesimal-main}), after introducing the
necessary background on measured laminations and grafting.  The
reduction to an infinitesimal statement is modeled on the argument of
Scannell-Wolf in \cite{scannell-wolf:grafting}, and uses Tanigawa's
properness theorem for grafting (Theorem \ref{thm:tanigawa}).

\sboldpoint{\S 3} describes shearing deformations of hyperbolic
surfaces, closely following the work of Bonahon on shearing
coordinates for Teichm\"uller space. Portions of this section are more
expository than is strictly necessary for the proofs of the of the
results, but we feel they help make the arguments easier to
understand.  The discussion culminates with the crucial
complex-linearity result of Bonahon (Theorem
\ref{thm:grafting-derivative}) that is used in the proof of
Theorem~\ref{thm:infinitesimal-main}.

\sboldpoint{\S 4} is devoted to the proofs of Theorems
\ref{thm:infinitesimal-main} and \ref{thm:main}, which follow using
the theory developed in \S\S 2-3.

\sboldpoint{\S 5} collects some applications of the main
theorem, including the proof of Corollary \ref{cor:graph} and a
rigidity result for quasi-Fuchsian manifolds.

\sboldpoint{\S 6} discusses the grafting coordinates for Teichm\"uller
space and the asymptotic behavior of extremal and hyperbolic length
functions on grafting rays; Theorems \ref{thm:extremal-length-intro}
and \ref{thm:hyp-length-intro} and associated derivative estimates are
proved here.

\subsection*{Acknowledgements.}
The authors thank Francis Bonahon, Howard Masur, Yair Minsky, and
Robert Penner for stimulating discussions related to this work.  They
also thank the referee for several suggestions that improved the
paper.

\section{Grafting and infinitesimal grafting}
\label{sec:grafting}

We begin with some background on measured laminations, grafting, and
tangentiability, which are needed to formulate the main technical
result (Theorem \ref{thm:infinitesimal-main}).

\subsection*{Laminations.}
As in the introduction, $S$ denotes a compact smooth surface of genus
$g >1$ and $\T(S)$ is the Teichm\"uller space of marked hyperbolic (or
conformal) structures on $S$.  We often use $X \in \T(S)$ to represent
a particular hyperbolic surface in a given marked equivalence class.

Let $\sS$ denote the set of free homotopy classes of simple closed
curves on $S$; we implicitly identify $\gamma \in \sS$ with its
geodesic representative on a hyperbolic surface $X \in \T(S)$.

A \emph{geodesic lamination} $\Lambda$ on a hyperbolic surface $X$ is
a foliation of a closed subset of $X$ by complete, simple geodesics.
Examples of geodesic laminations include simple closed hyperbolic
geodesics $\gamma \in \sS$ and disjoint unions thereof.

The notion of a geodesic lamination is actually independent of the
particular choice of $X$, in that a geodesic lamination on $X$
determines a geodesic lamination for any other hyperbolic structure $Y
\in \T(S)$ in a canonical way (see for example \cite[\S1]{bonahon:geodesic}).
Thus we speak of a geodesic lamination
on $S$, suppressing the choice of a particular metric.  Let $\GL(S)$
denote the set of all geodesic laminations on $S$ with the topology of
Hausdorff convergence of closed sets.

A geodesic lamination $\Lambda \in \GL(S)$ is \emph{maximal} if it is
not properly contained in another geodesic lamination, in which case
the complement of $\Lambda$ in $S$ is a union of ideal triangles.
Every geodesic lamination is contained in a maximal one, though not
necessarily uniquely.

\subsection*{Measured laminations.} A \emph{transverse measure} $\mu$ on
a geodesic lamination $\Lambda$ is an assignment of a positive Borel
measure to each compact transversal to $\Lambda$ in a manner
compatible with splitting and isotopy of transversals.  Such a measure
$\mu$ has \emph{full support} if there is no proper sublamination $\Lambda'
\subset \Lambda$ such that $\mu$ assigns the zero measure to
transversals disjoint from $\Lambda'$.

Let $\ML(S)$ denote the space of \emph{measured geodesic laminations}
on $S$, i.e. pairs $\lambda = (\Lambda, \mu)$ where $\Lambda \in
\GL(S)$ and $\mu$ is a transverse measure on $\Lambda$ of full
support.  We denote by $\lambda(\tau)$ the total measure assigned to a
transversal $\tau$ by $\lambda \in \ML(S)$.

The topology on $\ML(S)$ is that of weak-$*$ convergence of measures
on compact transversals.  The underlying geodesic lamination of
$\lambda \in \ML(S)$ is the \emph{support} of $\lambda$, written
$\supp(\lambda) \in \GL(S)$.  The space $\ML(S)$ has an action of
$\R^+$ by multiplication of transverse measures; the empty lamination
$0 \in \ML(S)$ is the unique fixed point of this action.

For any simple closed geodesic $\gamma \in \sS$, there is a measured
geodesic lamination (also $\gamma$) that assigns to a transversal
$\tau$ the counting measure on $\tau \cap \gamma$.  The rays $\{ t
\gamma \: | \: t \in \R^+, \: \gamma \in \sS\}$ determined by simple
closed curves are dense in $\ML(S)$.

The space $\ML(S)$ is a contractible topological manifold homeomorphic
to $\R^{6g-6}$, but it does not have a natural smooth structure.  Its
natural structure is that of a piecewise linear (PL) manifold, with
charts corresponding to train tracks.

Detailed discussion of the space $\ML(S)$ can be found in
\cite{thurston:notes} \cite{epstein-marden:earthquakes}
\cite{penner:train} \cite{otal:hyperbolization}.

\subsection*{Grafting.}
As mentioned in the introduction, Thurston showed that grafting along
simple closed curves has a natural extension to measured laminations,
giving a projective grafting homeomorphism $\Gr : \ML(S) \times \T(S)
\to \P(S)$ and a conformal grafting map $\gr : \ML(S) \times \T(S) \to
\T(S)$.  Tanigawa showed that the latter is a proper map when either
one of the two parameters is fixed:

\begin{thm}[{Tanigawa \cite{tanigawa:grafting}}]
\label{thm:tanigawa}
For any $\lambda \in \ML(S)$, the $\lambda$-grafting map
$\gr_\lambda : \T(S) \to \T(S)$ is proper.  For any $X \in \T(S)$, the 
$X$-grafting $\gr_{\param}X :
\ML(S) \to \T(S)$ is proper.
\end{thm}

The properness of the restricted grafting maps is used in the proofs of
Theorems \ref{thm:main} (in \S\ref{sec:proof}) and
\ref{thm:scannell-wolf} (in \cite{scannell-wolf:grafting}) to reduce a
global statement to a local one, which is then attacked using
infinitesimal methods.  In the case of $\lambda$-grafting, the
infinitesimal analysis is possible because $\gr_{\lambda} : \T(S) \to
\T(S)$ is differentiable, and even real-analytic
\cite[Cor.~2.11]{mcmullen:earthquakes}.  (A related real-analyticity
property along rays in $\ML(S)$ is discussed in \S6 below.)

The main step in \cite{scannell-wolf:grafting} is to show that the
differential map $d\!\gr_{\lambda} : T_X\T(S) \to T_{\gr_\lambda X}
\T(S)$ is an isomorphism.  Once that is established,
Theorem~\ref{thm:scannell-wolf} follows easily, since $\gr_\lambda$
is then a proper local diffeomorphism of $\T(S)$, hence a covering
map of the simply connected space $\T(S)$.

We will follow an analogous outline in the proof of
Theorem~\ref{thm:main}, but the infinitesimal analysis is complicated
by lack of smooth structure on $\ML(S)$, so the derivative of
$X$-grafting does not exist in the classical sense.  Instead we must
use a weaker notion of differentiability based on one-sided derivatives,
which we now discuss.

\subsection*{Tangentiability.}
A \emph{tangentiable} map\footnote{Some authors say instead that the
map is \emph{one-sided Gateaux differentiable}.} 
$f : U \to V$ between open sets in $\R^n$ is
a map with ``one-sided'' directional derivatives everywhere; in other
words, for each $x \in U$ and $v \in \R^n$, the limit
\begin{equation}
\dtzero f(x + t v) = \lim_{t \to 0^+} \frac{f(x+tv) - f(x)}{t}
\end{equation}
exists, and the convergence is locally uniform in $v$ (for equivalent
conditions, see \cite[\S 2]{bonahon:variations}).  This convergence
allows us to define $T_xf : \R^n \to \R^n$, the \emph{tangent map of f
at $x$}, by
\begin{equation*}
T_x f(v) = \dtzero f(x+tv)
\end{equation*}
Of course if $f$ is differentiable, then $T_x f$ is just the
derivative of $f$ at $x$, a linear map.  When $f$ is only
tangentiable, the map $T_x f$ is continuous and homogeneous in the
sense that $T_x f(\lambda v) = \lambda T_x f(v)$ for $\lambda \in
\R^+$ \cite[\S1]{bonahon:variations}.

A \emph{tangentiable manifold} is one whose transition functions are
tangentiable maps; examples include smooth manifolds and PL manifolds.
Thus $\ML(S)$ has a natural tangentiable structure.  The tangent space
$T_xM$ at a point $x$ of a tangentiable manifold is not naturally a
vector space, but has the structure of a cone.  The notion of
tangentiable map extends naturally to tangentiable manifolds.

We will say that a homeomorphism between two tangentiable manifolds is
a \emph{bitangentiable homeomorphism} if it and its inverse are
tangentiable, and if the tangent maps are everywhere homeomorphisms.
A convenient criterion for this is provided by the

\begin{lem}[{Bonahon \cite[Lem.~4]{bonahon:variations}}]
\label{lem:diffeo-criterion}
Let $f: M \to N$ be a homeomorphism between tangentiable manifolds.
If $f$ is tangentiable, and all of its tangent maps are injective,
then $f$ is a bitangentiable homeomorphism.
\end{lem}

Bonahon showed that grafting is compatible with that tangentiable
structure of $\ML(S)$ in the following sense:

\begin{thm}[{Bonahon \cite[Thm.~3]{bonahon:variations}}]
\label{thm:grafting-tangentiable}
The grafting map $\Gr:\ML(S) \times \T(S) \to \P(S)$ is a bitangentiable
homeomorphism.  In particular, the conformal grafting map $\gr:
\ML(S) \times \T(S) \to \T(S)$ is tangentiable, and for each $X \in
\T(S)$, the $X$-grafting map $\gr_{\param} X : \ML(S) \to \T(S)$ is
bitangentiable.
\end{thm}

In \cite{bonahon:variations}, Bonahon actually computes the tangent
map of grafting to show that grafting is tangentiable.  After
developing the shearing coordinates in \S3, Bonahon's description of
the tangent map (from which Theorem \ref{thm:grafting-tangentiable} is
derived) is given in Theorem \ref{thm:grafting-derivative}.

\subsection*{Differentiability.}
A curious feature of the tangentiability of grafting with respect to
$\ML(S)$ is that some fragments of classical differentiability remain.
For example, the inverse of the projective grafting map $\Gr^{-1} :
\P(S) \to \ML(S) \times \T(S)$ factors into the two projections
$p_\ML$ and $p_\T$ (as described in the introduction).  By Theorem
\ref{thm:grafting-tangentiable}, these are also tangentiable maps, but
since both the domain and range of $p_{\T}$ are smooth manifolds, it
makes sense to ask if this map is differentiable in the usual sense.
Extending Theorem \ref{thm:grafting-tangentiable}, Bonahon shows

\begin{thm}[{Bonahon \cite[Thm.~3]{bonahon:variations}}]
\label{thm:projection-differentiable}
The map $p_{\T} : \P(S) \to \T(S)$ is $C^1$.
\end{thm}

In the same article, Bonahon shows that $p_\T$ fails to be $C^2$ for a
certain family of punctured torus groups, which suggests that
$p_\T$ may fail to be $C^2$ for all  Teichm\"uller spaces $\T(S)$.

Finally we observe that $p_{\T}(\Gr_{\lambda} X) = X$, so for each
$\lambda \in \ML(S)$ we have $p_{\T} \circ \Gr_\lambda = \id$, that is,
the map $\Gr_{\lambda} : \T(S) \to \P(S)$ is a smooth section of
$p_{\T}$.  Since $\Gr : \ML(S) \times \T(S) \to \P(S)$ is a
homeomorphism, these sections fill up $\P(S)$, and we conclude

\begin{cor}
\label{cor:submersion}
The map $p_{\T}$ is a $C^1$ submersion.
\end{cor}

\subsection*{Infinitesimal $X$-grafting.}  Using the tangentiability of
grafting, we can formulate the infinitesimal statement that will be our
main tool in the proof of Theorem \ref{thm:main}:

\begin{thm}
\label{thm:infinitesimal-main}
The tangent map $T_\lambda \gr_{\param} X$ of the $X$-grafting map has
no kernel.  That is, if $\lambda_t$ is a tangentiable family of
measured laminations and $\dtzero \gr_{\lambda_t} X = 0$, then
$\dtzero \lambda_t = 0$.
\end{thm}

In the next two sections, we develop machinery to prove this result
about the derivative of grafting, then strengthen it to a local
injectivity result in order to prove Theorem~\ref{thm:main}.
Complications arise in both steps because the maps under consideration
are tangentiable rather than smooth.

\section{Shearing}
\label{sec:shearing}

In this section we describe the machinery of shearing cocycles for
geodesic laminations on a hyperbolic surface, borrowing heavily from
the papers of Bonahon \cite{bonahon:transverse},
\cite{bonahon:shearing}.  Some examples and discussion are included
here to clarify the technicalities that our later arguments will require.

\subsection*{Cocycles.}  Let $G$ be an abelian group and let $\Lambda \in
\GL(S)$.  A \emph{$G$-valued cocycle on $\Lambda$} is a map $\alpha$
that assigns to each transversal $\tau$ to $\Lambda$ an element
$\alpha(\tau) \in G$ in a manner compatible with splitting and
transversality-preserving isotopy.  The $G$-module of all $G$-valued
cocycles on $\Lambda$ is denoted $\sH(\Lambda,G)$.

Of particular interest for our purposes is the 
vector space of cocycles for
\emph{maximal} laminations with values in $\R$.  
While it is perhaps not clear from the definition, this vector space
is finite dimensional, and the dimension is the same for all maximal
laminations:

\begin{thm}[{Bonahon \cite[Prop.~1]{bonahon:shearing}}]
\label{thm:bonahon-cocycle}
Let $\Lambda \in \GL(S)$ be a maximal lamination.  Then
$\sH(\Lambda,\R) \simeq \R^{6g-6}$.
\end{thm}

The vector space $\sH(\Lambda,\R)$ carries a natural alternating
bilinear form $\omega : \sH(\Lambda,\R) \times \sH(\Lambda,\R) \to
\R$, the \emph{Thurston symplectic form}, which comes from the cup
product on $H^1(S, \R)$ (see \cite{penner:train} \cite[\S
3]{bonahon:shearing} \cite{bonahon-sozen:symplectic}).  When $\Lambda$
is maximal, the form $\omega$ is nondegenerate, making
$\sH(\Lambda,\R)$ a symplectic vector space.  While we do not use this
symplectic structure directly in the proof of the main theorem, it is
relevant to some of the constructions and examples in the sequel.

If $\lambda \in \ML(S)$ is a measured lamination with $\supp(\lambda)
\subset \Lambda \in \GL(S)$, then the total measure $\tau \mapsto
\lambda(\tau)$ defines a real-valued cocycle on $\Lambda$, which we
also denote by $\lambda \in \sH(\Lambda,\R)$.  Cocycles arising from
measures in this way take only nonnegative values on transversals;
Bonahon showed that the converse is also true.

\begin{thm}[{Bonahon \cite[Prop.~18]{bonahon:transverse}}]
\label{thm:positive-cocycle}
A transverse cocycle $\alpha \in \sH(\Lambda,\R)$ arises from a
transverse measure for $\Lambda$ if and only if $\alpha(\tau) \geq
0$ for every transversal $\tau$.
\end{thm}

We therefore define $\sM(\Lambda) \subset \sH(\Lambda,\R)$, the \emph{cone of
transverse measures} for $\Lambda$, to be the set of cocycles $\alpha$
satisfying $\alpha(\tau)\geq0$ for all transversals $\tau$.  The set
$\sM(\Lambda)$ is a convex cone in the vector space $\sH(\Lambda,\R)$.

While positive real-valued cocycles on $\Lambda$ correspond to
transverse measures, there is an essential difference between a
real-valued cocycle $\alpha \in \sH(\Lambda,\R)$ (whose value on a
transversal is a real number) and a \emph{signed transverse measure}
on $\Lambda$ which assigns a countably additive signed measure to each
transversal $\tau$ (see Examples \ref{ex:ray}-\ref{ex:spiral} below).

An analogous situation is the set function $[a,b] \mapsto (f(b) -
f(a))$, where $f$ is a real-valued function; this is a finitely
additive function on intervals, but it only arises from a signed Borel
measure if $f$ has bounded variation (which is automatic if $f$ is
monotone).  The connection between this example and real-valued
cocycles for a lamination can be seen through the ``distribution
function'' $f(x) = \alpha(\tau_x)$ where $\alpha \in \sH(\Lambda,\R)$,
$\tau : [0,1] \to S$ is a transversal, and $\tau_x =
\left.\tau\right|_{[0,x]}$.  This function is defined for a.e.~$x \in
[0,1]$ and monotonicity (for every $\tau$) is equivalent to $\alpha$ being
nonnegative.

The difference between measures and cocycles is also apparent from the
dimension of the span of $\sM(\Lambda)$, which is often $1$ (by the
solution to the Keane conjecture, see \cite{masur:ergodicity}
\cite{veech:ergodicity} \cite{rees:ergodicity}
\cite{kerckhoff:ergodicity}) and is never more than $(3g-3)$ (because
for maximal $\Lambda$, the space $\sH(\Lambda,\R) \simeq \R^{6g-6}$ is
a symplectic vector space in which $\Span(\sM(\Lambda))$ is isotropic,
see \cite{papadopoulos:dimension} \cite{levitt:dimension}
\cite{katok:dimension}).  It follows that $\sM(\Lambda)$ has positive
codimension when, for example, $\Lambda$ is maximal.

\subsection*{Shearing coordinates.}
Let $\Lambda \in \GL(S)$ be a maximal lamination, realized as a
partial foliation of $X \in \T(S)$ by hyperbolic geodesics.  Its lift
$\Tilde{\Lambda}$ to the universal cover $\Tilde{X} \simeq \H^2$
determines a (not necessarily locally finite) tiling of $\H^2$ by
ideal triangles.

A transversal $\tau : [a,b] \to \H^2$ to $\Tilde{\Lambda}$ determines
a pair of ideal triangles $T_a$ and $T_b$ which are the complementary
regions of $\Tilde{\Lambda}$ containing $\tau(a)$ and $\tau(b)$.  In
\cite{bonahon:shearing}, Bonahon constructs a \emph{shearing cocycle}
$\sigma(X) = \sigma^\Lambda(X) \in \sH(\Lambda, \R)$ from these data
with the property that $\sigma(X)(\tau)$ measures the ``relative
shear'' of the triangles $T_a$ and $T_b$ in $\H^2$.  For example, the
relative shear of two ideal triangles that share an edge is the signed
distance between the feet of the altitudes based on the common side.

Remarkably, $\sigma(X)$ determines the metric $X$, and the set of such
  cocycles admits an explicit description.  
Let $\sC(\Lambda) \subset \sH(\Lambda,\R)$ denote the set of
$\R$-valued cocycles that arise as shearing cocycles of hyperbolic
metrics, and recall that $\omega : \sH(\Lambda,\R) \times
\sH(\Lambda,\R) \to \R$ is the Thurston symplectic form.

\begin{thm}[{Bonahon \cite[Thm.~A,B]{bonahon:shearing}}]
\label{thm:shearing}
A cocycle $\alpha \in \sH(\Lambda,\R)$ is the shearing cocycle of a
hyperbolic metric if and only if only if $\omega(\alpha, \mu) > 0$ for
all $\mu \in \sM(\Lambda)$, and $\sC(\Lambda)$ is an open convex
cone with finitely many faces.  Furthermore $\sigma : \T(S) \to
\sC(\Lambda)$ is a real-analytic diffeomorphism.
\end{thm}

The condition $\omega(\alpha, \mu) > 0$ in Theorem \ref{thm:shearing}
is necessary because for each $\mu \in \sM(\Lambda)$, the Thurston
pairing $\omega(\sigma(X),\mu)$ is the hyperbolic length of $\mu$ on
$X$ \cite[Thm.~9]{bonahon:shearing}.

While the convex cone $\sC(\Lambda)$ has finitely many faces and is a
union of open rays in $\sH(\Lambda,\R)$, the zero cocycle
$0 \in \sH(\Lambda,\R)$ is \emph{not} an extreme point of $\sC(\Lambda)$. In
fact, the vector space of signed transverse measures (i.e.~$\Span
\sM(\Lambda)$) is $\omega$-isotropic (see
\cite{papadopoulos:dimension}), and therefore
\begin{equation}
\label{eqn:measures-boundary}
\Span \sM(\Lambda) \subset \partial \sC(\Lambda).
\end{equation}
A schematic representation $\sC(\Lambda)$ appears in Figure
 \ref{fig:shearing}, where $\partial \sC(\Lambda)$ contains a
 one-dimensional subspace of $\sH(\Lambda,\R)$.

\begin{remark} Thurston parametrizes $\T(S)$ by a convex
cone of shear coordinates in \S9 of \cite{thurston:stretch}.
The setting is in terms of duals to weights on train tracks,
but the results are equivalent to those in Theorem~\ref{thm:shearing}.
\end{remark}

\subsection*{Shearing maps.}
We now use the shearing embedding $\sigma : \T(S) \to \sH(\Lambda,\R)$
to turn translation in the vector space $\sH(\Lambda,\R)$ into a
(locally-defined) map of Teichm\"uller space.

Let $X \in \T(S)$ and $\alpha \in \sH(\Lambda,\R)$.  If the sum
$\sigma(X) + \alpha$ is the shearing cocycle of a hyperbolic surface,
we call this hyperbolic surface $\sh_\alpha X$, the \emph{shearing of $X$ by
$\alpha$}.  Thus $\sh_\alpha X$ is defined by the condition
\begin{equation}
\sigma(\sh_\alpha X) = \sigma(X) + \alpha.
\end{equation}
Since $\sigma(\T(S)) = \sC(\Lambda)$ is open, there is a neighborhood
$U \subset \sH(\Lambda,\R) \times \T(S)$ of $\{0\}\times\T(S)$ in
which the shearing map $\sh : U \to \T(S)$ is well-defined.

\begin{figure}
\begin{center}
\includegraphics[width=9cm]{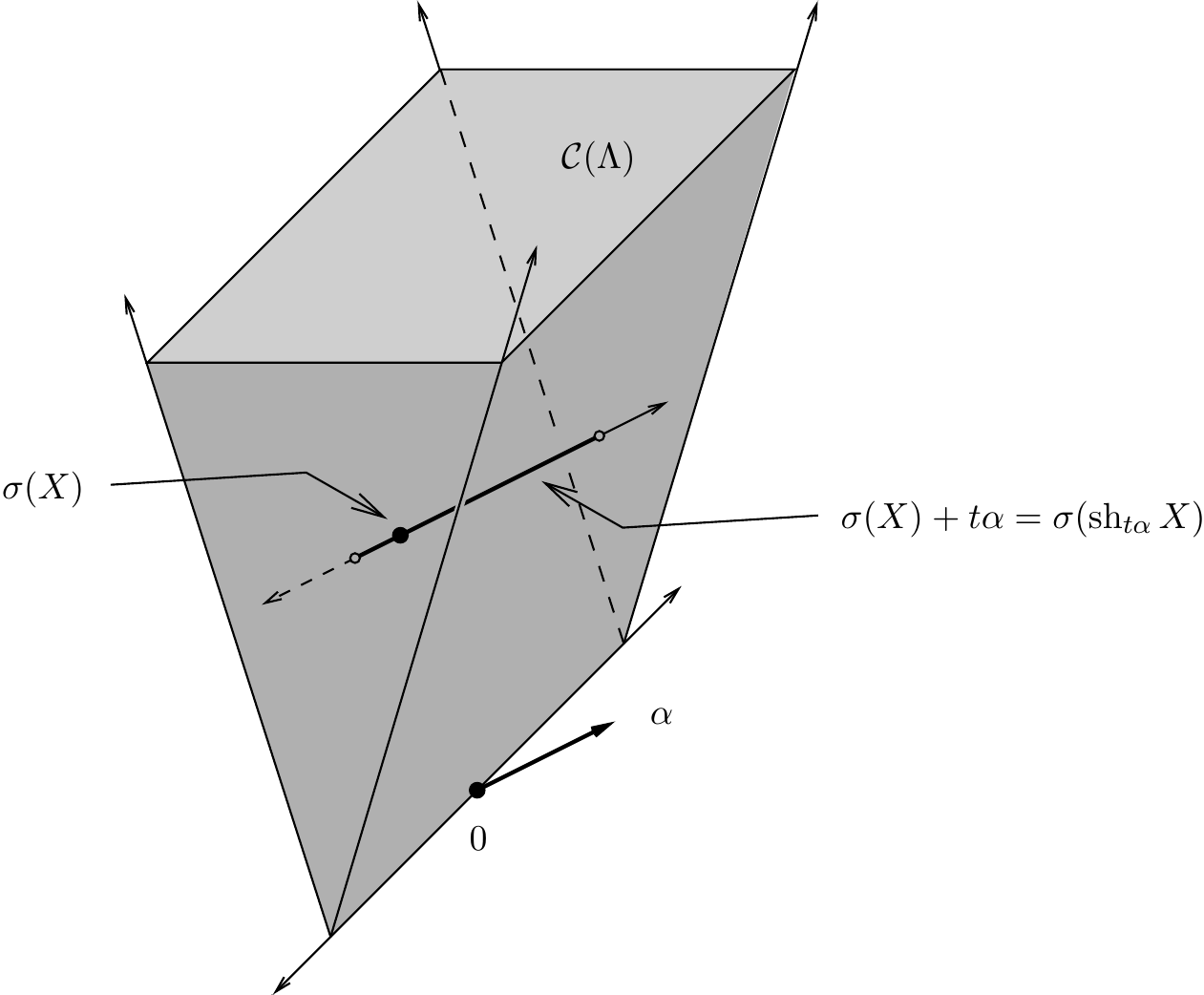}
\caption{The shearing map is a translation in the shearing embedding
  of $\T(S)$ in $\sH(\Lambda,\R)$.}
\label{fig:shearing}
\end{center}
\end{figure}

In particular, for any fixed $X \in \T(S)$ and $\alpha \in
\sH(\Lambda,\R)$ there is some $\epsilon > 0$ such that $\sh_{t
\alpha} X$ is defined for all $|t| < \epsilon$ (see Figure
\ref{fig:shearing}), and
\begin{equation*}
\dt \sigma \left(\sh_{t \alpha}X\right) = \alpha.
\end{equation*}

We consider a few examples of shearing maps to highlight the role of
$\sC(\Lambda)$ and the fact that $\sh_\alpha X$ is not defined for all
pairs $(\alpha,X)$.  First, if $\Lambda$ contains a simple closed
geodesic $\gamma$, then for each $t \in \R$ the cocycle $t \gamma$ is
a signed transverse measure.  Furthermore, for all $X \in \T(S)$, we
have $(\sigma(X) + t \gamma) \in \sC(\Lambda)$, since $t \gamma
\subset \sM(\Lambda)$ and the span of $\sM(\Lambda)$ is
$\omega$-isotropic; thus $\sh_{t \gamma} X$ is defined for all $t \in
\R$.  Concretely, the hyperbolic surface $\sh_{t \gamma}X$ is obtained
from $X$ by cutting along the geodesic $\gamma$ and then gluing the
two boundary components with a twist (by signed distance $t$).

This twisting example has a natural generalization: given a cocycle
$\lambda \in \sH(\Lambda,\R)$ representing a transverse measure for
$\Lambda$, the shearing $\sh_{t \lambda}(X)$ is again defined for all
$t \in \R$ and the resulting map $\sh_{t \lambda} : \T(S) \to \T(S)$
is called an \emph{earthquake}.  For further discussion of
earthquakes, see \cite{thurston:earthquakes}
\cite{kerckhoff:earthquakes} \cite{epstein-marden:earthquakes}
\cite{bonahon:earthquakes} \cite{mcmullen:earthquakes}.

As a final example, consider shearing a surface $X$ using a cocycle
$\alpha \in \sH(\Lambda,\R)$ that is itself the shearing cocycle of a
hyperbolic surface, i.e. $\alpha \in \sC(\Lambda)$.  For any
transverse measure $\lambda \in \sM(\Lambda)$, the hyperbolic length
of $\lambda$ on $\sh_{t \alpha} X$ can be computed using the Thurston
intersection form $\omega$ (cf.~\cite[Thm.~9]{bonahon:shearing}), and
we have
\begin{equation}
\ell(\lambda,\sh_{t \alpha} X) = \omega(\sigma(X) + t \alpha,\lambda) =
\ell(\lambda,X) + t\,\omega(\alpha,\lambda) = A + B t
\end{equation}
where $A,B >0$.  Since the length of a measured lamination is
positive, it follows that the set of $t$ for which $\sigma(X) + t
\alpha$ is the shearing cocycle of a hyperbolic metric (and hence
those for which $\sh_{t \alpha} X$ exists) is a subset of $\{ t >
-(A/B) \}$.

\subsection*{Tangent cocycles.}  Let  $\lambda_t \in \ML(S)$, $t \in [0,\epsilon)$ be a
tangentiable ray of measured laminations.  We will represent the
tangent vector $\dtzero \lambda_t$ by a transverse cocycle to a
certain geodesic lamination (as in \cite{bonahon:geodesic}).  We first
describe the underlying geodesic lamination.

The \emph{essential support of $\lambda_t$ at $t=0^+$} is a geodesic
lamination $\Lambda \in \GL(S)$ that reflects how the support of
$\lambda_t$ is changing for small positive values of $t$.  For a $PL$
family of measured laminations $\lambda_t$, the essential support is
the Hausdorff limit $\lim_{t \to 0^+} \supp(\lambda_t)$ of the
supporting geodesic laminations \cite[Prop.~4]{bonahon:geodesic}.  
For the general case, we only sketch the construction, and refer the
reader to \cite[\S2]{bonahon:geodesic} for details. 

First lift $\lambda_t$ to a family $\Tilde{\lambda_t}$ of measured
geodesic laminations in $\H^2$.  Define a set $\Tilde{\Lambda}$ of
geodesics in $\H^2$ as follows: a geodesic $\gamma$ belongs to
$\Tilde{\Lambda}$ if and only if for every smooth transversal $\tau :
[-\epsilon,\epsilon] \to \H^2$ with $\tau(0) \in \gamma$, the total
transverse measure of $\tau$ with respect to $\lambda_t$ is at least
$C t$ for some $C>0$ (depending on $\tau$ and $\gamma$) and all $t$
sufficiently small.  Then $\Tilde{\Lambda}$ is the lift of a geodesic
lamination $\Lambda \in \GL(S)$, the essential support of $\lambda_t$.

The tangent vector $\dtzero \lambda_t$ defines a real-valued
transverse cocycle $\dot{\lambda} \in \sH(\Lambda,\R)$ on the
essential support as follows:
\begin{equation*}
\dot{\lambda}(\tau) = \lim_{t \to 0^+} \frac{1}{t}
\left ( \lambda_t(\tau) - \lambda_0(\tau) \right ).
\end{equation*}
This is the \emph{tangent cocycle} of $\lambda_t$ at $t=0^+$.  Clearly
the same formula defines a cocycle for a lamination containing the
essential support of $\lambda_t$, so as a convenience we may assume
$\Lambda$ is maximal.  We illustrate this construction with two
examples:

\begin{example}
\label{ex:ray}
The cocycle determined by a measured lamination $\lambda \in
\sM(\Lambda)$ is the tangent cocycle of the ``ray'' $\lambda_t = (1 +
t)\lambda \in \ML(S)$ at $t=0^+$.
\end{example}

\begin{example}
\label{ex:spiral}
Let $S$ be a punctured torus with meridian $\alpha$ and longitude
$\beta$ (putting aside our assumption that $S$ is compact for a
moment), and consider the family $\lambda_t \in \ML(S)$ defined by the
conditions $\lambda_t(\alpha) = 1$ and $\lambda_t(\beta) = t$.  Thus
for each $n \in \N$, the measured lamination $\lambda_{1/n}$ is a
simple closed curve with homology class $(n[\alpha] + [\beta])$ and
weight $1/n$.

The essential support of $\lambda_t$ at $t=0^+$
is the geodesic lamination $\Lambda = \alpha \cup \eta$, where $\eta$
is an infinite simple geodesic that spirals toward $\alpha$ in each
direction.  The derivative $\dot{\lambda} \in \sH(\Lambda,\R)$ is a
cocycle with full support and indefinite sign (compare Bonahon's
example \cite[p.~104]{bonahon:geodesic}).
\end{example}

This second example shows that the tangent cocycle may not be a
transverse measure, and the essential support may not admit a measure
of full support.  This illustrates some of the difficulties of using
differential methods on the space of measured laminations and of
adapting the methods in \cite{scannell-wolf:grafting} to prove the
main theorem.

Following Thurston, Bonahon showed that the association of a cocycle
$\dot{\lambda}$ to a tangentiable family $\lambda_t$ provides a linear
model for each $PL$ face of the tangent space $T_{\lambda_0} \ML(S)$
(see \cite{bonahon:geodesic}, \cite[\S6]{thurston:stretch}).  When
$\supp(\lambda_0)$ is not maximal, however, no single maximal
lamination $\Lambda$ can be chosen to contain the essential support of
\emph{every} family $\lambda_t$, and so there is no embedding
$T_{\lambda_0} \ML(S) \to \sH(\Lambda,\R)$.

In contrast, when $\supp(\lambda_0)=\Lambda$ is maximal (a generic
situation that excludes, for example, closed leaves), the tangent
cocycle construction defines a homeomorphism $T_{\lambda_0}\ML(S)
\simeq \sH(\Lambda,\R)$, giving the tangent space a canonical linear
structure.  This was observed in \cite{thurston:stretch}.

\subsection*{Complex linearity.} So far we have seen real-valued
cocycles on geodesic laminations arise in two different contexts:
first as shearing cocycles providing coordinates for $\T(S)$, and then
as tangent vectors to families of measured laminations.  Building on
these two constructions, the following result of Bonahon will allow us
to connect the derivative of grafting with respect to $\T(S)$ and
$\ML(S)$:

\begin{thm}[{Bonahon \cite[Prop. 5]{bonahon:variations} \cite[\S10]{bonahon:shearing}}]
\label{thm:grafting-derivative}
Let $X \in \T(S)$ and $\lambda \in \ML(S)$.  For each maximal geodesic
lamination $\Lambda \in \GL(S)$ containing the support of $\lambda$,
there is a complex-linear map $L = L(\Lambda,\lambda,X) :
\sH(\Lambda,\C) \to T_{\Gr_{\lambda}X} \P(S)$ which determines the
tangent map of $\Gr$ in tangent directions carried by $\Lambda$, in
the following sense:

Let $\lambda_t \in \ML(S)$ be a tangentiable family of measured
laminations with $\lambda_0 = \lambda$ and with essential support
contained in $\Lambda$, and let $\dot{\lambda} = \dtzero \lambda_t \in
\sH(\Lambda,\R)$.  Let $X_t \in \T(S)$ be a smooth family of
hyperbolic structures with $X_0 = X$ whose derivative in the $\Lambda$
shearing embedding is $\dot{\sigma} = \dtzero \sigma(X_t)$.

Then $t \mapsto \Gr_{\lambda_t} X_t$ is a tangentiable curve in
$\P(S)$ and
\begin{equation*}
\dtzero \Gr_{\lambda_t} = L(\dot{\sigma} + i \dot{\lambda})
\end{equation*}
Similarly $\gr: \ML(S) \times \T(S) \to \T(S)$ is tangentiable and
its derivative has the same complex linearity property.
\end{thm}

In terms of the piecewise linear structure of $\ML(S)$, Theorem
\ref{thm:grafting-derivative} says that on each linear face of the
tangent space $T_{(\lambda,X)}(\ML(S) \times \T(S))$, the tangent map
$T_{(\lambda,X)} \gr$ is the restriction of a complex-linear map
$\sH(\Lambda,\C) \to T_{\gr_\lambda X} \T(S)$ (see
\cite[\S2]{bonahon:variations} \cite[\S10]{bonahon:shearing}).

\section{Proof of the main theorem}
\label{sec:proof}

With the necessary background in place, we can now show that the
tangent map of $\gr_{\param}X$ has no kernel.

\begin{proof}[Proof of Theorem \ref{thm:infinitesimal-main}.]
Fix $X$ and suppose that $\dtzero \gr_{\lambda_t} X = 0$.  Let
$\dot{\lambda} \in \sH(\Lambda,\R)$ denote the derivative of
$\lambda_t$ at $t=0^+$, with $\Lambda \in \GL(S)$ maximal.  We must
show that $\dot{\lambda}=0$.

For all $t$ sufficiently small, the shearing $X_t = \sh_{t
\dot{\lambda}} X \in \T(S)$ is defined and satisfies $X_0 = X$ and
$\dot{\lambda} = \dt \sigma(X_t)$.  Therefore by Theorem
\ref{thm:grafting-derivative},
\begin{equation*}
i \left [ \dtzero \gr_{\lambda_0} X_t \right ] = i L(\dot{\lambda}) =
  L(i \dot{\lambda}) = \dtzero \gr_{\lambda_t} X_0 = 0,
\end{equation*}
where $L$ is the complex-linear map representing the tangent map of
$\gr$ in tangent directions carried by $\Lambda$.

By Theorem \ref{thm:scannell-wolf}, $\gr_{\lambda_0}$ is an immersion,
and so $\dtzero X_t = 0$.  Thus
\begin{equation*}
\dot{\lambda} = \dt \sigma(X_t)
= 0.
\end{equation*}
\end{proof}

If the grafting map were continuously differentiable in the usual
sense, the proof of the main theorem would now be straightforward,
using linearity of the derivative and the inverse function theorem
to conclude that $\gr_{\param}X$ is a local diffeomorphism.  We will
follow this general outline, but we will need to use additional
properties of the grafting maps to strengthen the infinitesimal result
to a local one.

In fact, some argument specific to grafting is necessary at this
point.  In general, a tangent map that has no kernel need not be
injective (consider $\C \to \R^+$ by $z \mapsto |z|$ at $z=0$), and
even if the tangent map is injective, tangentiability does not imply
the continuous variation of derivatives needed for the inverse
function theorem.

With these potential problems in mind, we analyze the $X$-grafting map
$\gr_{\param}X$ as the composition of the projective $X$-grafting map
$\Gr_{\param} X$ and the smooth projection $\pi : \P(S) \to \T(S)$.
Let $N = \dim \T(S) = \dim \ML(S) = \frac{1}{2} \dim \P(S) = 6g-6$.

Let $G_X = \{ \Gr_\lambda X \: | \: \lambda \in \ML(S) \} \subset
\P(S)$ denote the image of $\Gr_{\param} X$.  Recall the map $p_\T :
\P(S) \to \T(S)$ is defined by $p_\T(\Gr_\lambda X) = X$, and so $G_X
= p_{\T}^{-1}(X)$ is a fiber of this map.  Since $p_{\T}$ is a $C^1$
submersion (by Corollary \ref{cor:submersion}), the set $G_X$ is
actually a $C^1$ submanifold of $\P(S)$ of dimension $N = \dim \P(S) -
\dim \T(S)$.  In particular, $G_X$ is smoother than its tangentiable
parameterization by $\ML(S)$ would suggest.

We can recast Theorem \ref{thm:infinitesimal-main} as a result about
the tangent space to $G_X$ as follows:

\begin{thm}
\label{thm:transverse}
For any $X \in \T(S)$, the $C^1$ submanifold $G_X \subset \P(S)$ is
transverse to the map $\pi : \P(S) \to \T(S)$; that is, for any $Z \in
G_X$, we have $T_Z G_X \cap (\ker d\pi) = \{ 0 \}$.  Equivalently, the
distributions $\ker d\pi$ and $\ker dp_{\T}$ in $T \P(S)$ are
transverse.
\end{thm}

\begin{proof}
Suppose not, i.e.~that there exists $v \in (T_Z G_X \cap \ker d\pi)$
with $v \neq 0$.  Since $Z \in G_X$, we have $Z = \Gr_{\lambda} X$ for
some $\lambda \in \ML(S)$.  Choose a $C^1$ path $Z_t$ in $G_X$ with
$Z_0 = Z$ and $\dtzero Z_t = v$, and let $\Gr^{-1}(Z_t) = (\lambda_t,
X)$.  Since $\Gr^{-1}$ is a bitangentiable homeomorphism, the
family $\lambda_t$ is tangentiable and satisfies $\dtzero \lambda_t
\neq 0$.  On the other hand, $\dtzero \gr_{\lambda_t}X = \dtzero
\pi(Z_t) = d \pi(\dtzero Z_t) = 0$, contradicting Theorem
\ref{thm:infinitesimal-main}.
\end{proof}

While a general tangent map can have no kernel and yet fail to be
injective, the fact that $G_X$ is $C^1$ (and thus has linear tangent
spaces) rules out this behavior for $\Gr_{\param}X$:

\begin{thm}
\label{thm:tangent-homeomorphism}
The tangent map $T_\lambda \gr_{\param} X$ of the conformal $X$-grafting map is
a homeomorphism.
\end{thm}

\begin{proof}
We study the diagram of grafting maps:
\begin{equation*}
\xymatrix@R+3mm@C+8mm{
\ML(S) \ar[r]^-{\Gr_{\param}X} \ar[rd]_{\gr_{\param}X} &
**[r] G_X \subset \P(S)
\ar[d]^{\pi} \\
& \T(S)
}
\end{equation*}

Fix $\lambda \in \ML(S)$ and for brevity let $f = T_\lambda
\gr_{\param} X$ be the tangent map of the conformal $X$-grafting map;
our goal is to show that $f$ is a homeomorphism.

Similarly, let $F = T_\lambda \Gr_{\param} X$ be the tangent map of the
projective $X$-grafting map.  Then $f = d \pi \circ F$, where $d\pi$
is the differential of $\pi : \P(S) \to \T(S)$ at $Z =
\Gr_{\lambda}X$, i.e. the tangent maps form a corresponding diagram:
\begin{equation*}
\xymatrix@R+3mm@C+8mm{
T_{\lambda} \ML(S) \ar[r]^-{F}
\ar[rd]_{f} &
**[r] T_{Z}G_X \subset T_{Z} \P(S)
\ar[d]^{d \pi} \\
& T_{\pi(Z)} \T(S)
}
\end{equation*}
Note that $F$ and $f$ are continuous homogeneous maps,
while $d \pi$ is linear.

Since $\Gr$ is a tangentiable homeomorphism, its restriction
$\Gr_{\param}X: \ML(S) \to \P(S)$ is a tangentiable injection,
i.e.~its tangent map $F$ is both injective and a homogeneous
homeomorphism onto its image.  Since $G_X$ is a $C^1$ submanifold of
$\P(S)$, the image of $F$ is the linear subspace $T_{Z} G_X \subset
T_{Z} \P(S)$.

First we show that $f$ injective.  Suppose on the contrary there exist
$\dot{\lambda}_1, \dot{\lambda}_2 \in T_\lambda \ML(S)$, distinct and
nonzero, and $f(\dot{\lambda}_1) = f(\dot{\lambda}_2)$.  Then $v =
F(\dot{\lambda}_1) - F(\dot{\lambda}_2) \in T_Z G_X$ is nonzero since
$F$ is injective, and so $v \in \ker d \pi$, which contradicts Theorem
\ref{thm:transverse}.

Thus $f : T_\lambda \ML(S) \to T_{\gr_\lambda X} \T(S)$ is injective.
Since $f$ is also a homogeneous map 
between cones of the
same dimension, it is a homeomorphism.
\end{proof}

\begin{remark}
The fact that $T_\lambda \Gr_{\param} X$ is a homeomorphism onto a
linear subspace of $T_{\Gr_\lambda X} \P(S)$ is also observed by
Bonahon in the proof of Prop.~12 in \cite{bonahon:variations}.
\end{remark}

Finally we complete our study of $\Gr_{\param}X$ and $\gr_{\param}X$
by proving the main theorem:

\begin{proof}[Proof of Theorem \ref{thm:main}.]
Let us consider the restriction of the forgetful map $\pi : \P(S) \to
\T(S)$ to the $C^1$ submanifold $G_X \subset \P(S)$.  Since $G_X$ is
$N$-dimensional, this projection is a local diffeomorphism at $Z \in
G_X$ if the subspaces $T_Z G_X$ and $\ker d \pi$ of $\T_Z \P(S)$ are
transverse (by the inverse function theorem).  By Theorem
\ref{thm:transverse}, this is true for every $Z \in G_X$, so
$\left.\pi\right|_{G_X}$ is a local $C^1$ diffeomorphism.

Thus the conformal $X$-grafting map is the composition of the
homeomorphism $\Gr_{\param} X : \ML(S) \to G_X$ and the local
homeomorphism $\left.\pi\right|_{G_X} : G_X \to \T(S)$, so
$\gr_{\param} X$ is a local homeomorphism.  As we noted in Section
\ref{sec:grafting}, by the properness of $\gr_{\param}X$ (Theorem
\ref{thm:tanigawa}), it follows that this map is a homeomorphism.
Thus $X$-grafting is a tangentiable homeomorphism with injective
tangent maps (Theorem \ref{thm:tangent-homeomorphism}), which by Lemma
\ref{lem:diffeo-criterion} is a bitangentiable homeomorphism.
\end{proof}

\section{Applications}
\label{sec:applications} 

In this section, we collect some applications of the main theorem
itself and of the techniques used in its proof, and we discuss some
related questions about grafting coordinates and $\CP^1$ structures.

\subsection*{Projections.}
We begin by proving the main corollary of Theorems \ref{thm:main} and
\ref{thm:scannell-wolf} about the grafting coordinates for a fiber
$P(X) = \pi^{-1}(X) \subset \P(S)$:

\begin{main_corollary}
For each $X\in\T(S)$, the space $P(X)$ of $\CP^1$ structures on $X$ is
a graph over each factor in the grafting coordinate system.  In fact,
we have:
\begin{enumerate}
\item The projection $\left.p_{\T}\right|_{P(X)} : P(X) \to \T(S)$ is a $C^1$
diffeomorphism.
\item The projection $\left.p_{\ML}\right|_{P(X)}:P(X)\to\ML(S)$ is a
bitangentiable homeomorphism.
\end{enumerate}
\end{main_corollary}

\begin{proof}
First we consider the regularity of the maps.  Bonahon showed that
$p_{\T}$ is $C^1$ (Theorem \ref{thm:projection-differentiable}), and
$p_{\ML}$ is tangentiable because it is the composition of $\Gr^{-1}$,
a bitangentiable homeomorphism (Theorem
\ref{thm:grafting-tangentiable}), and the projection to one factor of
a product of tangentiable manifolds.
\begin{enumerate}
\item \label{item:pt} \emph{Proof that $\left.p_{\T}\right|_{P(X)}$ is a
diffeomorphism.} First of all, the map $\left.p_{\T}\right|_{P(X)}$
is a homeomorphism because it has inverse map
\begin{equation*}
\label{eqn:teich-projection-inverse}
Y \mapsto \Gr_{(\gr_{\param}X)^{-1}(Y)}Y
\end{equation*}
where the map $(\gr_{\param}X)^{-1} : \T(S) \to \ML(S)$ exists by
Theorem \ref{thm:main}.

Thus it suffices to show that $\left.p_{\T}\right|_{P(X)}$ is a local
diffeomorphism.  But the kernel of $\left.d p_{\T}\right|_{P(X)}$ is the
intersection of $\ker d \pi$ and $\ker d p_{\T}$ in $T \P(S)$, which
is zero by Theorem \ref{thm:transverse}.  So the derivative of
$\left.p_{\T}\right|_{P(X)}$ is an isomorphism, and by the inverse
function theorem this map is a local diffeomorphism.

\item \label{item:pml} \emph{Proof that $\left.p_{\ML}\right|_{P(X)}$
  is a bitangentiable homeomorphism.}  As in \eqref{item:pt}, we first
  show that $\left.p_{\ML}\right|_{P(X)}$ is a homeomorphism by
  exhibiting an inverse map, 
\begin{equation*}
\label{eqn:ml-projection-inverse}
\lambda \mapsto \Gr_{\lambda}(\gr_{\lambda}^{-1}(Y))
\end{equation*}
where $\gr_\lambda^{-1}$ exists by Theorem \ref{thm:scannell-wolf}.  

By Lemma \ref{lem:diffeo-criterion}, we need only show that the
tangent map of $\left.p_{\ML}\right|_{P(X)}$ is everywhere injective.
That this tangent map has no kernel also follows easily from Theorem
\ref{thm:scannell-wolf}, but to show injectivity we will use an
argument modeled on the proofs of Theorems
\ref{thm:infinitesimal-main} and \ref{thm:tangent-homeomorphism}.

Suppose on the contrary that two distinct, nonzero tangent vectors
$v_1,v_2 \in T_ZP(X)$ have the same image in $T_\lambda \ML(S)$ under
the tangent map of $p_\ML$, where $Z = \Gr_\lambda Y$.  Then
differentiable paths in $P(X)$ with tangent vectors $v_1$ and $v_2$
are mapped by $\Gr^{-1}$ to tangentiable paths $\Gr_{\lambda_t^{(1)}}
Y_t^{(1)}$ and $\Gr_{\lambda_t^{(2)}} Y_t^{(2)}$, respectively, where
$Y^{(k)}_0 = Y$ and $\lambda^{(k)}_0 = \lambda$ for $k=1,2$.

In this notation, the image of $v_k$ by the tangent map $T_Z p_\ML$ is
$\dtzero \lambda^{(k)}_t$, so we have $\dtzero \lambda^{(1)}_t =
\dtzero \lambda^{(2)}_t = \dot{\lambda}$.  In particular there is a
single geodesic lamination $\Lambda$ containing the essential support
of both families $\lambda^{(k)}_t$, and their common tangent vector
defines a cocycle $\dot{\lambda} \in \sH(\Lambda,\R)$.  Using the
shearing embedding of Teichm\"uller space in $\sH(\Lambda,\R)$ gives
cocycles $\dot{\sigma}_k = \dtzero(\sigma^\Lambda(Y_t^{(k)}))$, and
since $v_1 \neq v_2$ we have $\dot{\sigma}_1 \neq \dot{\sigma}_2$.

By Theorem \ref{thm:grafting-derivative}, there is a complex-linear
map $L : \sH(\Lambda,\C) \to T_X \T(S)$ that gives the tangent map of
$\Gr$ for tangent vectors to $\ML(S) \times \T(S)$ at $(\lambda,X)$
representable by complex-valued cocycles on $\Lambda$, so $v_k =
L(\dot{\sigma}_k + i \dot{\lambda})$.  Now consider $\dot{\sigma} =
\dot{\sigma}_1 - \dot{\sigma}_2 \neq 0 $, which is the tangent vector
to the shearing family $Y_t = \sh_{t \dot{\sigma}} Y$; we have
\begin{equation*}
\begin{split}
\dtzero \Gr_\lambda Y_t &= L(\dot{\sigma})\\
& = L\left ( (\dot{\sigma}_1 + i \dot{\lambda}) - (\dot{\sigma}_2 + i
\dot{\lambda}) \right ) \\
& = v_1 - v_2 \in T_Z P(X),
\end{split}
\end{equation*}
and so $\dtzero \gr_\lambda Y_t = 0$, which by Theorem
\ref{thm:scannell-wolf} implies that $\dtzero Y_t = 0$, and so
$\dot{\sigma} = \dtzero \sigma^\Lambda(Y_t) = 0$, a contradiction.
Thus the tangent map of $p_\ML$ is injective, as required.
\end{enumerate}
\end{proof}

\begin{remark}
Since Bonahon constructs an explicit example to show that the full
projection map $p_\T : \P(S) \to \T(S)$ is not necessarily $C^2$, it
would be interesting to know if a similar construction could be used
to show that $\left . p_\T \right |_{P(X)}$ need not be $C^2$.
\end{remark}

\subsection*{Pruning.} 
Since the $\lambda$-grafting map $\gr_\lambda : \T(S) \to \T(S)$ is a
homeomorphism (Theorem \ref{thm:scannell-wolf}), there is an inverse
map $\pr_\lambda : \T(S) \to \T(S)$ which we call \emph{pruning by
  $\lambda$}.  Roughly speaking, grafting by $\lambda$ inserts a Euclidean
subsurface along the leaves of $\lambda$, and pruning by $\lambda$
removes it.  Allowing $\lambda$ to vary, we obtain the pruning map
$\pr : \ML(S) \times \T(S) \to \T(S)$, and fixing $X$ we have the
$X$-pruning map $\pr_{\param} X : \ML(S) \times \{X\} \to \T(S)$.

We can reformulate Corollary \ref{cor:graph} in terms of pruning as follows:

\begin{cor}
For each $X \in \T(S)$, the $X$-pruning map
$\pr_{\param}X:\ML(S)\to\T(S)$ is a bitangentiable homeomorphism.
\end{cor}

\begin{proof}
The graph of the $X$-pruning map consists of the pairs $(\lambda,Y)$
such that $\gr_\lambda Y = X$, which is simply the fiber $P(X) \subset \P(S)
\simeq \ML(S) \times \T(S)$.  In terms of the projections projections
$p_\ML : P(X) \to \ML(S)$ and $p_\T : P(X) \to \T(S)$, we have
$\pr_{\param} X = p_\T \circ p_\ML^{-1}$, which is a bitangentiable
homeomorphism by Corollary \ref{cor:graph}.
\end{proof}

Previously it was known that $\pr_{\param} X$ is a ``rough
homeomorphism'', i.e. a proper map of degree $1$, and that is has a
natural extension to the Thurston compactification of $\T(S)$ and the
projective compactification of $\ML(S)$ by $\PML(S) = (\ML(S) -
\{0\})/\R^+$.  Furthermore, the resulting boundary map $\PML(S) \to
\PML(S)$ is the \emph{antipodal map} relative to $X$.  For details,
see \cite{dumas:antipodal}.

\subsection*{Internal Coordinates for the Bers slice} 
\begin{figure}
\begin{center}
\includegraphics[width=10cm]{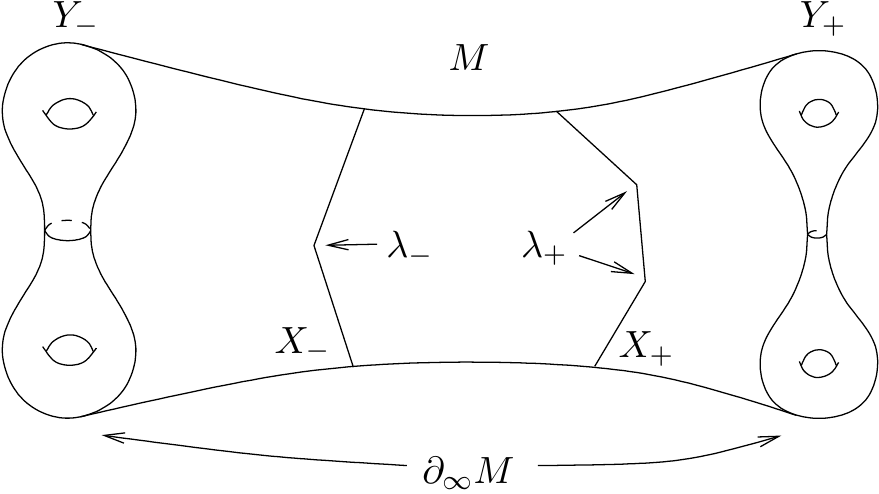}
\caption{Geometric data associated to a quasi-Fuchsian manifold $M \in
  \QF$.  In a Bers slice, $Y_-$ is fixed.}
\label{fig:qf}
\end{center}
\end{figure}

Let $\QF = \QF(S)$ denote the space of marked quasi-Fuchsian
hyperbolic 3-manifolds homeomorphic to $S \times \R$ (see
\cite{bers:simultaneous-uniformization} \cite{nag:teichmuller}).  Each
such manifold $M \in \QF$ has ideal boundary $\partial_\infty M = Y_+
\sqcup Y_-$ with conformal structures $Y_{\pm} \in \T(S)$ and convex
core boundary surfaces $X_{\pm} \in \T(S)$ with bending laminations
$\lambda_\pm \in \ML(S)$ (represented schematically in Figure
\ref{fig:qf}).  Furthermore, the convex core and ideal boundary
surfaces satisfy $\gr_{\lambda_{\pm}} X_{\pm} = Y_{\pm}$
(\cite[Thm.~2.8]{mcmullen:earthquakes}, see also
\cite{kamishima-tan:grafting}).

In his celebrated holomorphic embedding \cite{bers:embedding} of the
Teichm\"uller space $\T(S)$ into $\C^{3g-3}$, Bers focused on the
slice $\B_Y\subset\QF$ of those quasi-Fuchsian manifolds with a fixed
ideal boundary surface $Y_- = Y$.  Bers showed that a quasi-Fuchsian
manifold $M \in \B_Y$ is uniquely determined by its other conformal
boundary surface $Y_+ \in \T(S)$, which can be chosen arbitrarily (the
\emph{Simultaneous Uniformization Theorem}
\cite{bers:simultaneous-uniformization}).

A corollary of our main theorem is that $\B_Y$ may also be
parameterized by the hyperbolic structure $X_-$:

\begin{cor} 
\label{cor:qf}
Let $M,M' \in \QF$ be marked quasi-Fuchsian manifolds, and suppose an
end of $M$ (respectively $M'$) has ideal boundary $Y$ (resp. $Y'$) and
the associated convex core boundary surface has hyperbolic metric $X$
(resp. $X'$).  If $Y$ and $Y'$ are conformally equivalent and $X$ and
$X'$ are isometric, then $M$ is isometric to $M'$.
\end{cor}

\begin{proof}
By hypothesis $Y'=Y$ and $X'=X$ as points in Teichm\"uller space, so
the bending measures $\lambda$ and $\lambda'$ of the convex core
boundaries satisfy $\gr_\lambda X =\gr_{\lambda'}X=Y$.  By Theorem
\ref{thm:main}, we have $\lambda=\lambda'$.

A quasi-Fuchsian manifold $M$ is uniquely determined up to isometry by
the hyperbolic metric $X$ and bending lamination $\lambda$ of one of
its convex core boundary surfaces, since one can use $X$ and $\lambda$
to construct the associated equivariant pleated plane in $\H^3$ and
its holonomy group $\pi_1 M \subset \PSL_2(\C)$ (see
\cite{epstein-marden:earthquakes}).  As
$M$ and $M'$ share these data, they are isometric.
\end{proof}

For another perspective on this corollary, we can consider the Bers
slice $\B_Y$ as a subset of $P(Y)$, where a quasi-Fuchsian manifold $M
\in \B_Y$ is identified with the projective structure it induces on
its ideal boundary surface $Y$.  Then Corollary \ref{cor:qf} is equivalent to
the statement that for every $Y \in \T(S)$, the projection map $\left
. p_\T \right |_{\B_Y} : \B_Y \to \T(S)$ is injective.  Of course this
can also be derived from Corollary \ref{cor:graph}, which shows,
furthermore, that $\left . p_T \right |_{\B_Y}$ is a $C^1$ embedding
(since $\B_Y \subset P(Y)$ is open). 

\begin{remarks}
\item In \cite{scannell-wolf:grafting}, 
it was observed that a manifold $M \in \B_Y$ is also determined by
the bending lamination $\lambda_- \in \ML(S)$ on the same side as the
fixed conformal structure $Y$.  It is not known whether either the
bending lamination $\lambda_+$ or the hyperbolic structure $X_+$
determine elements of $\B_Y$.

\item More generally, one can ask what geometric data determine $M \in
\QF$ up to isometry.  Bonahon and Otal \cite{bonahon-otal:bending}
showed that if $\lambda_+$ and $\lambda_-$ bind the surface and are
supported on simple closed curves, then the pair
$(\lambda_+,\lambda_-)$ determines $M \in \QF$ uniquely.  Recently,
Bonahon \cite{bonahon:almost-fuchsian} showed that this restriction to
simple closed curves may be lifted for elements of $\QF$ which are
sufficiently close to the Fuchsian subspace $\F\subset\QF$, and Series
\cite{Se04} proved that $\lambda_+$ and $\lambda_-$ determine $M$ when
$S$ is a once-punctured torus.
\end{remarks}

\section{Grafting coordinates and rays}
\label{sec:rays}

In this final section we discuss how the ray structure of $\ML(S)$ is
transported to $\T(S)$ by the $X$-grafting map.

Recall (from Section \ref{sec:grafting}) that the action of $\R^+$ on
$\ML(S)$ by scaling transverse measures gives this space the structure
of a cone, with the empty lamination $0$ as its base point.  By
Theorem \ref{thm:main}, for each $X \in \T(S)$ we can use the
$X$-grafting map to parameterize the Teichm\"uller space by $\ML(S)$,
providing a global system of ``polar coordinates'' centered at $X$.
In this coordinate system, the ray $\R^+ \lambda \in \ML(S)$
corresponds to the \emph{grafting ray} $\{ \gr_{t \lambda} X \: | \: t
\in \R^+\}$, a properly embedded path starting at $X$, and
Teichm\"uller space is the union of these rays.

It would be interesting to understand the geometry of this coordinate
system, and especially the grafting rays.  Thus we ask:
\begin{question}
What is the behavior of the grafting ray $t \mapsto \gr_{t \gamma} X$
and how does it depend on $\lambda$ and $X$?
\end{question}
Theorems \ref{thm:extremal-length} and \ref{thm:hyp-length} below
address this question by estimating the extremal and
hyperbolic lengths of the grafting lamination (and their derivatives)
along a ray.
The asymptotic behavior of certain grafting rays in
relation to the Teichm\"uller metric has also been investigated by
D\'{i}az and Kim, see \cite{diaz-kim}.
 
Naturally, one first wonders about the regularity of grafting rays,
since the grafting map itself exhibits a combination of tangentiable
and differentiable behavior.  However, along rays the grafting map is
as smooth as possible:

\begin{thm}[{McMullen \cite{mcmullen:earthquakes}}]
\label{thm:rays-analytic}
The grafting ray $t \mapsto \gr_{t \lambda}$ is a real-analytic map
from $\R^+$ to $\T(S)$; in fact, it is the restriction of a
real-analytic map $\{ t \geq -\epsilon \} \to \T(S)$ for some
$\epsilon>0$ (depending on $\lambda$).  Furthermore, if $\lambda_n \to
\lambda$ in $\ML(S)$, then the $\lambda_n$ grafting rays converge in
$C^\omega$ to the $\lambda$ grafting ray.
\end{thm}

\begin{remarks}
\item In \cite{mcmullen:earthquakes} it is shown that for any $\lambda \in
\ML(S)$, the \emph{complex earthquake map} $\mathrm{eq}_\lambda : \H
\to \T(S)$ is holomorphic and extends to an open neighborhood of $\H
\cup \R$ (Proposition 2.6 and Theorem 2.10).  Furthermore these maps
vary continuously with $\lambda \in \ML(S)$ (Theorem 2.5), as do their
derivatives, since they are holomorphic.  Since the grafting ray is
the restriction of the complex earthquake to $i \R$, Theorem
\ref{thm:rays-analytic} follows immediately.

\item The regularity of grafting rays and complex earthquakes is closely
related to (and in part, an application of) the analyticity of
quake-bend deformations of surface group representations in
$\PSL_2(\C)$ established by Epstein and Marden
\cite{epstein-marden:earthquakes}.
\end{remarks}

\subsection*{Extremal length.}
When $\lambda$ is supported on a simple closed geodesic $\gamma$, it
is clear that a large grafting will result in a surface in which
$\gamma$ has small extremal length, because $\gr_{t \gamma} X$
contains an annulus of large modulus homotopic to $\gamma$.  Refining
this intuition, we will establish the estimate:

\begin{thm}
\label{thm:extremal-length}
For each $X \in \T(S)$, the extremal length of $\lambda$ is of order
$1/t$ on the $\lambda$-grafting ray $t \mapsto \gr_{t \lambda}X$
and is monotone decreasing for all $t \gg 0$.  Specifically, we have
\begin{enumerate}
\item $E(\lambda, \gr_{t \lambda}X) =  \frac{\ell(\lambda,X)}{t} +
  O(t^{-2})$, and  \label{item:bound}
\item $\frac{d}{dt} E(\lambda, \gr_{t \lambda} X) =
\frac{-\ell(\lambda,X)}{t^2} + O(t^{-3})$
\label{item:monotone}
\end{enumerate}
where $E(\lambda,Y)$ denotes the extremal length of $\lambda$ on the
Riemann surface $Y$.  The implicit constants depend only on $\chi(S)$.
\end{thm}

Note that Theorem \ref{thm:extremal-length} includes the results
stated in the introduction as Theorem
\ref{thm:extremal-length-intro}.  Before giving the proof of
Theorem \ref{thm:extremal-length}, we fix notation and recall some
concepts from Teichm\"uller theory used therein.

\subsection*{Annuli.}
Let $A$ be an annular Riemann surface of modulus $M$ and let $E(A) =
1/M = E(\gamma,A)$ be the extremal length of $\gamma$, the nontrivial
isotopy class of simple closed curves on $A$.  Then $A$ is isomorphic
to a rectangle $R_A = [0,E(A)] \times (0,1) \subset \C$ with its
vertical sides identified.  We call the complex local coordinate $z$
on $A$ coming from this realization the \emph{natural coordinate} for
$A$.  Similarly, the induced flat metric $|dz|$ on $A$ is the
\emph{natural metric}, with respect to which $A$ is a Euclidean
cylinder of height $1$ and circumference $E(A)$.

\subsection*{Jenkins-Strebel differentials.}
For any isotopy class $\gamma$ of simple closed curves on a compact
Riemann surface $X$, there is a unique embedded annulus $A_\gamma
\subset X$ homotopic to $\gamma$ of maximum modulus $\Mod(A_\gamma) =
1/E(\gamma,X)$.  The annulus $A_\gamma$ is dense in $Y$, and if $z$ is
the natural coordinate for $A_\gamma$, the quadratic differential
$dz^2$ on $A$ extends holomorphically to a quadratic differential on
$X$, the \emph{Jenkins-Strebel differential} for $\gamma$.

\subsection*{Foliations.}
A holomorphic quadratic differential $\phi$ on a Riemann surface $X$
has an associated singular \emph{horizontal foliation} $\F(\phi)$
whose leaves integrate the distribution of tangent vectors $v$
satisfying $\phi(v) \geq 0$.  Integration of $|\Im \sqrt \phi|$ gives
a transverse measure on $\F(\phi)$.  Similarly $\F(-\phi)$ is the
\emph{vertical foliation}, whose transverse measure comes from $|\Re
\sqrt \phi|$.

When $\phi$ is a Jenkins-Strebel differential on a compact surface,
the nonsingular leaves of $\F(\phi)$ are closed and homotopic to
$\gamma$; in the realization of the Jenkins-Strebel annulus $A_\gamma$
as a rectangle with identifications, these are the horizontal lines,
while leaves of $\F(-\phi)$ are vertical lines.  The transverse
measures for $\F(\phi)$ and $\F(-\phi)$ are given by $|dy|$ and $|dx|$
in the rectangle, respectively.
Thus the closed leaves of $\F(\phi)$ have total measure $E(A_\gamma) =
E(\gamma,X)$ with respect to the transverse measure of the vertical
foliation $\F(-\phi)$.  Furthermore these closed leaves realize the
minimum transverse measure among all curves homotopic to $\gamma$ (see
\cite[Lem.~11.5.3]{gardiner:teichmuller-theory}).

\subsection*{Pairing and extremal length.}
The natural pairing between Beltrami differentials $\mu =
\mu(z) \slashdzbardz$ (with $\mu(z) \in L^\infty$) and integrable holomorphic
quadratic differentials $\phi = \phi(z) dz^2$ on a Riemann surface
$X$ is given by
\begin{equation*}
\langle \mu, \phi \rangle = \Re \int_X \mu \phi = \int_X \mu(z)
\phi(z) \:|dz|^2.
\end{equation*}
When $X=A$ is an annulus with natural coordinate $z_A$, pairing a
Beltrami differential with $\phi = dz_A^2$ gives the infinitesimal
change in extremal length $E(A)$ (see
\cite[\S1.9]{gardiner:teichmuller-theory}); that is, if $A_t$ is a
family of annuli identified by a family of quasiconformal maps with
derivative $\mu_t$, then
\begin{equation}
\label{eqn:extremal-length-derivative}
\frac{d}{dt}E(A_t) = 2 \langle \mu(t), dz_{A_t}^2 \rangle.
\end{equation}
Similarly, when $\phi$ is a Jenkins-Strebel differential on a compact
surface $X$, pairing with $\phi$ gives the differential of the
extremal length function on Teichm\"uller space:
\begin{thm}[{Gardiner \cite[Thm.~11.8.5]{gardiner:teichmuller-theory}}]
\label{thm:gardiner}
Let $t \mapsto X_t \in \T(S)$ be a differentiable path whose tangent
vector is represented by the Beltrami differential $\mu(t)$ on $X_t$.
Let $\gamma$ be an isotopy class of simple closed curves and $E(t) =
E(\gamma,X_t)$ its extremal length on $X_t$.  Then
\begin{equation*}
E'(t) = 2 \langle \mu(t), \phi(t) \rangle
\end{equation*}
where $\phi(t)$ is the Jenkins-Strebel differential for $\gamma$ on
$X_t$.
\end{thm}

\subsection*{Stretching annuli.}
Let $X$ be a compact Riemann surface and $A \subset X$ an annulus in
the homotopy class of $\gamma$, a simple closed curve.  The natural
coordinate $z$ on $A$ gives a Beltrami differential $\slashdzbardz$ on $A$,
which extends to a Beltrami differential on $X$ by setting it to zero
on $(X-A)$.  This differential represents an infinitesimal affine
stretch of $A$.  We will be interested in the extent to which
$\slashdzbardz$ affects the extremal length of $\gamma$, as estimated in

\begin{lem}
\label{lem:stretching-annulus}
Let $A \subset X$, $\gamma$, and $z$ be as above, and let $\phi$ be
the Jenkins-Strebel differential on $X$ for $\gamma$.  Then
\begin{equation*}
E(X) \geq \left \langle \dzbardz, \phi \right \rangle \geq \frac{2 E(X)^2}{E(A)} - E(X)
\end{equation*}
where $E(A) = E(\gamma,A)$ and $E(X) = E(\gamma,X)$ are the extremal
lengths of $\gamma$ on these two surfaces.
\end{lem}

\begin{remarks}
\item When $A = A_\gamma$ is the Jenkins-Strebel annulus (of maximum
modulus), we have $E(A) = E(X)$ and both inequalities in Lemma
\ref{lem:stretching-annulus} become equalities.  However this is clear
since $\slashdzbardz = \overline{\phi}/|\phi|$ if $z$ is the natural
coordinate of the Jenkins-Strebel annulus.  The point of the Lemma is that we
also have tight bounds for the pairing when $A$ has \emph{nearly}
maximum modulus.

\item A related estimate for nearly maximal annuli is used in
Kerckhoff's proof of that that foliation map $\F : Q(X) \to \ML(S)$ is
a homeomorphism, see
\cite[Lem.~3.2]{kerckhoff:asymptotic}.
\end{remarks}

\begin{proof}
Throughout the proof we use the natural coordinate $z$ to identify $A$
with a rectangle $R_A \subset \C$ whose vertical sides are identified.

Writing the restriction of $\phi$ to $A$ in terms of the natural
coordinate, we have
\begin{equation*}
\phi = \phi(z) dz^2 = (\alpha + i \beta)^2 dz^2,
\end{equation*}
where $\alpha$ and $\beta$ are real-valued functions defined locally up to a
common sign away from the zeros of $\phi$; in particular the functions 
$\alpha^2$, $\beta^2$ and $|\alpha|$ are well-defined almost
everywhere.  We want to estimate
\begin{equation*}
\begin{split}
\left \langle \dzbardz, \phi \right \rangle &= \Re \int_A \dzbardz \: \phi(z) dz^2\\
&=
\int_{A} \Re \left ( (\alpha + i \beta)^2\right ) \:
|dz|^2\\
& = \int_0^1
\int_0^{E(A)} (\alpha^2 - \beta^2) dx dy.
\end{split}
\end{equation*}
Note that the pairing is computed as an integral over $A$, rather
than $X$, because the Beltrami coefficient $\slashdzbardz$ is understood to
be zero on $(X-A)$.

We first derive the upper bound on the pairing. Since $A$ is a subset
of $X$, we have
\begin{equation}
\label{eqn:alphasq-betasq}
\int_0^1 \int_0^{E(A)} (\alpha^2 + \beta^2) dx dy = \int_A |\phi| \leq
\int_X |\phi| = E(X).
\end{equation}
Since $\beta^2 \geq 0$, the same upper bound applies to the integral
of $\alpha^2 - \beta^2$, giving
\begin{equation*}
\left \langle \dzbardz, \phi \right \rangle = \int_0^1 \int_0^{E(A)} (\alpha^2 - \beta^2) dx
dy \leq E(X).
\end{equation*}

To establish the lower bound on the pairing, note that the horizontal
lines in $R_A$ represent closed curves in $X$ homotopic to $\gamma$,
so the total transverse measure of any one of these with respect to
$\F(-\phi)$ is at least $E(X)$.  The transverse measure of a curve is
the integral of $|\Re \sqrt{\phi}| = |\alpha dx - \beta dy|$, but
$dy=0$ on horizontal lines, so we have
\begin{equation*}
\int_0^{E(A)} |\alpha(x + i y)| dx \geq E(X)
\end{equation*}
for all $y \in [0,1]$.  Integrating over $y$ and applying the
Cauchy-Schwarz inequality, we obtain
\begin{equation}
\label{eqn:alphasq}
\int_0^1 \int_0^{E(A)} \alpha^2 dx dy \geq \frac{E(X)^2}{E(A)}.
\end{equation}

Multiplying \eqref{eqn:alphasq} by $2$ and subtracting \eqref{eqn:alphasq-betasq} we have
\begin{equation*}
\langle \dzbardz, \phi \rangle = \int_0^1 \int_0^{E(A)} (\alpha^2 - \beta^2) dx
dy \geq \frac{2 E(X)^2}{E(A)} - E(X).
\end{equation*}
\end{proof}
Note that the proof of Lemma \ref{lem:stretching-annulus} is
essentially a calculation on the annulus $A$ and uses little about the
enclosing surface $X$ except that it is foliated by closed
trajectories of $\phi$.  Indeed, the same argument can be applied with
the compact surface $X$ replaced by an annulus $B$ and with $\phi =
dz_B^2$ the natural quadratic differential on $B$, and it is this
version we will need in the proof of Theorem \ref{thm:hyp-length}:
\begin{lem}
\label{lem:stretching-subannulus}
Let $B$ be an annular Riemann surface of finite modulus and $A \subset
B$ a homotopically essential subannulus.  Then
\begin{equation*}
E(B) \geq \left \langle \dbard{z_A}, dz_B^2 \right \rangle \geq \frac{2
  E(B)^2}{E(A)} - E(B)
\end{equation*}
\end{lem}

\subsection*{Extremal length and grafting rays.}
Using Lemma \ref{lem:stretching-annulus} as the main technical tool,
we are now ready to give the

\begin{proof}[Proof of Theorem \ref{thm:extremal-length}]
  First, we consider the case when $\lambda = \gamma$ is a simple
  closed geodesic with unit weight.  For brevity let $Y_t = \gr_{t
    \gamma} X$; we abbreviate $E(t) = E(\gamma,Y_t)$ and $\ell =
  \ell(\gamma,X)$.

The proof of \eqref{item:bound} follows the usual pattern for an
extremal length estimate (see \cite[\S3]{kerckhoff:asymptotic}): a
particular annulus homotopic to $\gamma$ bounds $E(t)$ from above,
while a particular conformal metric on the surface bounds $E(t)$ from
below.  In this case the annulus is the grafting cylinder $A_t \subset
Y_t$ of modulus $t/\ell$, and the conformal metric on $Y_t$ is the
\emph{Thurston metric}--the union of the product metric on $A_t =
[0,t] \times \gamma$ and the hyperbolic metric of $X$ (see
\cite[\S2.1]{tanigawa:grafting}).  Applying the geometric and analytic
definitions of extremal length gives
\begin{equation*}
\frac{\ell}{t} > E(\gamma,Y_t) > \frac{\ell^2}{t \ell + A}
> \frac{\ell}{t} - \frac{A}{t^2}
\end{equation*}
where $A = 4 \pi (g-1)$ and $(t \ell + A)$ is the area of the Thurston
metric; thus \eqref{item:bound} follows.

For statement \eqref{item:monotone}, we must estimate the derivative of
extremal length.  For all $s,t > 0$, there is a natural quasiconformal
map from $Y_t$ to $Y_{s}$ that is affine on $A_t$, stretching it
vertically in the natural coordinate, and conformal on $(Y_t - A_t)$.
The derivative of this family of maps at $s=t$ is the Beltrami
differential
\begin{equation*}
\mu(t) = -(2t)^{-1} \dbard{z_t} 
\end{equation*}
where $z_t$ is the natural coordinate on the grafting annulus $A_t$
and the Beltrami differential $\slashdbard{z_t}$ is understood to be
identically zero outside $A_t$.

By Theorem \ref{thm:gardiner}, the derivative of extremal length along
the grafting ray is
\begin{equation}
\label{eqn:derivative-formula}
\frac{d}{dt} E(t) = 2 \langle \mu(t), \phi(t)
\rangle = - t^{-1} \left \langle \dbard{z_t}, \phi(t) \right \rangle.
\end{equation}
We estimate the pairing $\langle \slashdbard{z_t}, \phi(t) \rangle$ using
Lemma \ref{lem:stretching-annulus}; starting with the upper bound, we
have
\begin{equation}
\label{eqn:pairing-upper}
\left \langle \dbard{z_t}, \phi(t) \right \rangle \leq E(t) \leq
\frac{\ell}{t}
\end{equation}
while the lower bound from the lemma gives
\begin{equation*}
\begin{split}
\left \langle \dbard{z_t}, \phi(t) \right \rangle \geq \frac{2
  E(t)^2}{E(A_t)} - E(t) = \frac{2 t}{\ell} E(t)^2 - E(t).
\end{split}
\end{equation*}
Using $E(t) \geq \ell/t - A/t^2$ on the first term and $E(t) \leq
\ell/t$ on the second, we obtain
\begin{equation}
\label{eqn:pairing-lower}
\left \langle \dbard{z_t}, \phi(t) \right \rangle 
\:\geq\:
\frac{2 t}{\ell}\left (\frac{\ell}{t} - \frac{A}{t^2} \right )^2 -
  \frac{\ell}{t}
\:\geq\:
\frac{\ell}{t} - \frac{4 A}{t^2}.
\end{equation}
Multiplying \eqref{eqn:pairing-upper} and \eqref{eqn:pairing-lower} by
$-1/t$ and using the formula \eqref{eqn:derivative-formula} for the
derivative of extremal length, we have
\begin{equation*}
-\frac{\ell(\gamma,X)}{t^2} \leq \frac{d}{dt} E(\gamma, \gr_{t
  \gamma} X) \leq -\frac{\ell(\gamma,X)}{t^2} + \frac{4 A}{t^3},
\end{equation*}
which gives part \eqref{item:monotone}, completing the proof of Theorem
\ref{thm:extremal-length} for simple closed geodesics.

The limiting arguments that extend (1) and (2) to general laminations
are completely analogous, so we will only give details for the
former.  Given a measured lamination $\lambda \in \ML(S)$, define
$$ \delta(s, \lambda, X) = s^2 \left | E(\lambda, \gr_{s \lambda} X ) -
  \frac{\ell(\lambda,\gr_{s \lambda}X)}{s} \right |.$$ Since
hyperbolic length, extremal length, and grafting are continuous on
$\ML(S) \times \T(S)$, this is a continuous nonnegative function
$\delta : \R^+ \times \ML(S) \times \T(S) \to \R$.  Part (1) of the
Theorem is equivalent to the statement that $\delta$ is bounded.  We
have shown that this is true for simple closed geodesics, i.e. there
exists $C > 0$ depending only on $\chi(S)$ such that for all
simple closed curves $\gamma$ and all $s>0$, we have
$$ \delta(s,\gamma, X) \leq C.$$
Suppose that $\lambda = c \gamma$ is a weighted simple closed
geodesic.  Since hyperbolic length scales linearly in the transverse
measure, while extremal length scales quadratically, we have
\begin{equation*}
\begin{split}
\delta(s, c \gamma, X) &= s^2 \left | E(c \gamma, \gr_{s c \gamma} X ) -
\frac{\ell(c \gamma,\gr_{s c \gamma}X)}{s} \right |\\
&= s^2 c^2 \left | E(\gamma, \gr_{s c \gamma} X) -
  \frac{\ell(\gamma,\gr_{s c \gamma}X)}{s c} \right |\\
&= \delta(s c, \gamma, X) \leq C.
\end{split}
\end{equation*}
Since weighted simple closed geodesics are dense in $\ML(S)$, and
$\delta$ is continuous, this shows that $\delta \leq C$, establishing
(1) for all measured laminations.

A similar limiting argument applies to part (2), where one bounds the
error function
$$ \Delta(s, \lambda, X) = s^3 \left | \left (\frac{d}{ds}
E(\lambda,\gr_{s \lambda} X) \right ) + \frac{\ell(\lambda,X)}{s^2}
\right |.$$ In this case, the function $\Delta$ in continuous because
extremal length is $C^1$ on Teichm\"uller space (see
\cite{gardiner-masur:extremal-length}), with its derivative varying continuously in
$\ML(S)$, and since the derivatives of grafting rays depend
continuously on the lamination (by Theorem \ref{thm:rays-analytic}).
\end{proof}

As a consequence of Theorem \ref{thm:extremal-length}, extremal length
decreases along grafting rays outside of a compact set in $\T(S)$ of
the form $\{ \gr_\lambda X \: | \: \ell(\lambda,X) \leq L \}$.  

\subsection*{Hyperbolic length.} 
The hyperbolic length of $\lambda$ along a grafting ray is more
difficult to control than the extremal length, but for the case of a
single curve with large weight, the same techniques used in the proof
of Theorem \ref{thm:extremal-length} give:

\begin{thm}
\label{thm:hyp-length}
Let $X \in \T(S)$ and let $\gamma$ be a simple closed hyperbolic
geodesic on $X$.  Then the hyperbolic length of $\gamma$ is of order
$1/t$ on the $\gamma$-grafting ray and is monotone decreasing for all
$t \gg 0$.  Specifically, we have
\begin{enumerate}
\item $\ell(\gamma, \gr_{t \gamma} X) = \frac{ \pi \ell(\gamma,
X)}{t} + O(t^{-2})$ \label{item:hyp-bound}
\item $\frac{d}{dt} \ell(\gamma, \gr_{t \gamma} X) = -\frac{\pi
  \ell(\gamma,X)}{t^2} + O(t^{-3})$ \label{item:hyp-monotone}
\end{enumerate}
as $t \to \infty$, where the implicit constants depend on $X$ and
$\gamma$.
\end{thm}

Note that Theorem \ref{thm:hyp-length} includes the results stated in
the introduction as Theorem \ref{thm:hyp-length-intro}.

\begin{proof}
As before let $Y_t = \gr_{t \gamma} X$, and abbreviate $\ell(t) =
\ell(\gamma,Y_t)$ and $E(t) = E(\gamma,Y_t)$.

A standard argument using the collar lemma shows that the extremal
and hyperbolic length of a curve are asymptotically proportional when
the length is small; specifically, we have
\begin{equation}
\pi E(t) - C E(t)^2 < \ell(t) < \pi E(t)
\end{equation}
for all $t$ such that $E(t) < 1$, where $C$ is a universal constant.
Since $E(t) < 1$ for all $t > \ell(\gamma,X)$, substituting the
estimate for $E(t)$ from Theorem \ref{thm:extremal-length} gives
\eqref{item:hyp-bound}.

Now we establish the derivative estimate \eqref{item:hyp-monotone}.
Let $\hat{Y}_t$ denote the cover of $Y_t$ corresponding to the
subgroup $\langle \gamma \rangle \subset \pi_1(Y_t)$.  Since
$\hat{Y}_t$ is conformally equivalent to an annulus of modulus $\pi /
\ell(t)$, we have $\ell'(t) = \pi \dt E(\hat{Y}_t)$.

Recall that $\mu(t) = (-2t)^{-1} \slashdbard{z_t}$ represents the
derivative of $Y_t$, where $z_t$ is the natural coordinate on $A_t$.
Thus the derivative of the annular covers $\hat{Y}_t$ is represented
by the
pullback Beltrami differential $p^*(\mu(t))$ where $p : \hat{Y}_t \to
Y_t$ is the covering projection.  Let $w_t$ denote the natural
coordinate on $\hat{Y}_t$.  By \eqref{eqn:extremal-length-derivative}
we have
\begin{equation}
\label{eqn:lprime-is-pairing}
\ell'(t) = \pi \frac{d}{dt} E(\hat{Y}_t) = 2 \pi \left \langle p^*(\mu(t)), dw_t^2 \right \rangle.
\end{equation}

To estimate this pairing, we analyze the differential $p^*(\mu(t))$;
its support is the preimage of the grafting cylinder $p^{-1}(A_t)
\subset \hat{Y}_t$, which consists of
\begin{rmenumerate}
\item A homotopically essentially annulus $\hat{A}_t$ such
  that $\left.p\right|_{\hat{A}_t} : \hat{A}_t \to A_t$ is a conformal
  isomorphism, and
\item A complementary set $\Omega = p^{-1}(A_t) - \hat{A_t}$ that is a
 disjoint union of countably many simply connected regions $\Sigma_i
 \subset \hat{Y}_t$ such that the restriction of $p$ to any one of
 them gives a universal covering $\left.p\right|_{\Sigma_i} : \Sigma_i
 \to A_t$,
\end{rmenumerate}
as depicted in Figure \ref{fig:cover}.  Therefore we have
\begin{equation}
\label{eqn:pairing-terms}
\langle p^*(\mu(t)), dw_t^2 \rangle = 
\int_{\hat{A}_t} p^*(\mu(t)) dw_t^2 + \int_{\Omega} p^*(\mu(t)) dw_t^2,
\end{equation}
and we can analyze these two terms individually.

\begin{figure}
\begin{center}
\includegraphics{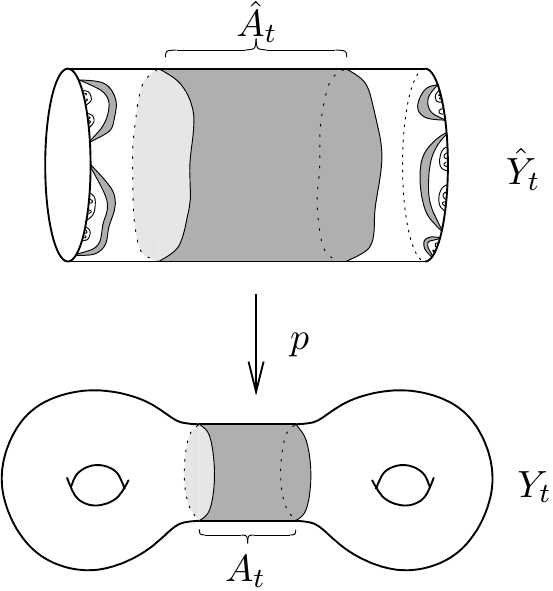}
\caption{The grafting annulus and its lifts to the annular cover.}
\label{fig:cover}
\end{center}
\end{figure}

By the length estimates of part \eqref{item:hyp-bound}, the reciprocal
moduli $E(\hat{A}_t) = E(A_t) = \ell(0)/t$ and $E(\hat{Y}_t) = \ell(t)
/ \pi$ differ by $O(t^{-2})$.  It follows that the subannulus
$\hat{A}_t$ accounts for nearly all of the area of $\hat{Y}_t$.  In
fact, restricting the metric $|dw_t|$ to $\hat{A}_t$ and using the
analytic definition of extremal length gives
\begin{equation*}
\area(\hat{A}_t, |dw_t|) \geq
\frac{\ell(\gamma,|dw_t|)^2}{E(\hat{A}_t)} =
\frac{E(\hat{Y}_t)^2}{E(\hat{A}_t)} \geq \frac{\ell(0)}{t} -
\frac{C}{t^2},
\end{equation*}
where $C$ depends on $\gamma$ and $X$.
Since $\area(\hat{Y}_t,|dw_t|) = \ell(0)/t$ and $\Omega$ is disjoint
from $\hat{A}_t$, we have $\area(\Omega, |dw_t|) \leq C/t^2$.

This area estimate implies that the second term in
\eqref{eqn:pairing-terms} is negligible, i.e.
\begin{equation}
\label{hyp-pairing-abs}
\begin{split}
\left |
 \int_{\Omega} p^*(\mu(t)) dw_t^2 \right | &\leq \|p^*(\mu(t))\|_\infty  \area(\Omega, |dw_t|)\\
& \leq \left ( \frac{1}{2t} \right ) \left ( \frac{C_2}{t^2} \right ) = O(t^{-3}).
\end{split}
\end{equation}

Now we consider the first term in \eqref{eqn:pairing-terms}.  Since
$\left.p\right|_{\hat{A}_t}$ is a conformal isomorphism, we have
$\left.p^*(\mu(t))\right|_{\hat{A}_t} = (-2t)^{-1} \slashdbard{\hat{z}_t}$
where $\hat{z}_t$ is the natural coordinate of $\hat{A}_t$.  Applying
Lemma \ref{lem:stretching-subannulus} we have
\begin{equation*}
E(\hat{Y}_t) \geq \left \langle \dbard{\hat{z}_t}, dw_t^2 \right \rangle \geq
  \frac{2 E(\hat{Y}_t)^2}{E(\hat{A}_t)} - E(\hat{Y}_t).
\end{equation*}
As before we substitute $E(\hat{Y}_t) = \ell(t) / \pi$,
$E(\hat{A}_t) = \ell(0) / t$, and apply the estimates for $\ell(t)$ to
obtain
\begin{equation*}
\int_{\hat{A}_t} p^*(\mu(t)) dw_t^2 = (-2t)^{-1} \left \langle
\dbard{\hat{z}_t}, dw_t^2 \right \rangle = - \frac{\ell(0)}{2t^2} +
O(t^{-3}).
\end{equation*}

Thus we have estimates for both terms in \eqref{eqn:pairing-terms},
and applying the formula \eqref{eqn:lprime-is-pairing} for $\ell'(t)$
gives the desired result:
\begin{equation*}
\ell'(t) = - \frac{\pi \ell(0)}{t^2} + O(t^{-3}).
\end{equation*}
\end{proof}

\begin{remark}
In \cite[Cor.~3.2]{mcmullen:earthquakes}, McMullen shows that
$\ell(\lambda,\gr_{t \lambda} X) < \ell(\lambda,X)$ for all $X \in
\T(S)$, $\lambda \in \ML(S)$ and $t>0$; it is also mentioned that
this upper bound can be strengthened to
\begin{equation*}
\ell(\gamma,\gr_{t \gamma}X) \leq \frac{\pi}{\pi+t} \ell(\gamma,X),
\end{equation*}
for any simple closed curve $\gamma$.  For details on this and a
corresponding lower bound, see \cite[Prop.~3.4]{diaz-kim}.  Theorem
\ref{thm:hyp-length} shows that this upper bound is asymptotically
sharp.
\end{remark}

\end{document}